\documentclass[11pt]{article}

\usepackage{amsfonts}
\usepackage{epsfig,graphics,subfigure,psfrag,amsmath,amssymb}
\usepackage{latexsym,wrapfig,picinpar,picins}
\usepackage{bm}
\usepackage{amsthm}
\usepackage{graphicx}
\usepackage{mathrsfs}
\usepackage{fancyref}
\usepackage{dcolumn}
\usepackage{color}
\usepackage{pstricks}
\usepackage{epstopdf}
\usepackage{subfigure}
\usepackage{multicol}
\usepackage{enumerate}
\usepackage{caption}
\usepackage{geometry}
\usepackage{threeparttable}
\makeatletter
\@addtoreset{equation}{section}
\makeatother
\renewcommand{\theequation}{\arabic{section}.\arabic{equation}}

\newtheorem{thm}{Theorem}[section]

\newtheorem{lem}[thm]{Lemma}

\newtheorem{rem}[thm]{Remark}
\def\E\mathbb{ E}

\newcommand{\be}{\begin{eqnarray}}
\newcommand{\ee}{\end{eqnarray}}
\newcommand{\ben}{\begin{eqnarray*}}
\newcommand{\een}{\end{eqnarray*}}

\begin{document}

\title{ Numerical methods for  Porous Medium Equation by an Energetic Variational Approach}

\author{Chenghua Duan\footnotemark[2]\and Chun Liu\footnotemark[3]\and Cheng Wang\footnotemark[4]\and Xingye Yue\footnotemark[1]}

 \renewcommand{\thefootnote}{\fnsymbol{footnote}}
 \footnotetext[2]{Department of Mathematics, Soochow University, Suzhou 215006, China (chduan@stu.suda.edu.cn).}
 \footnotetext[3]{Department of Applied Mathematics, Illinois Institute of Technology, Chicago, IL 60616, USA,
 (cliu124@iit.edu).}

 \footnotetext[4]{Department of Mathematics, University of Massachusetts, Dartmouth, North Dartmouth, MA, 02747-2300, USA (cwang1@umassd.edu).}
 \footnotetext[1]{Corresponding author. Department of Mathematics, Soochow University, Suzhou 215006, China (xyyue@suda.edu.cn).}
%



\date{}
\maketitle

{
\begin{abstract}
\footnotesize
We study numerical methods for porous media equation (PME).  There are two important characteristics: the finite  speed propagation of the free boundary and the potential waiting time, which make the problem not easy to handle. Based on different dissipative energy laws, we develop two numerical schemes by an energetic variational approach. Firstly, based on  $f \log f$ as the total energy form of the dissipative law,   we obtain the trajectory equation, and then construct a fully discrete scheme. It is proved that the scheme is  uniquely solvable on an admissible  convex set by taking the advantage of the singularity of the total energy.   Next, based on $\frac{1}{2 f}$ as  the total energy form of the dissipation law, we construct a linear numerical  scheme for the corresponding trajectory equation. Both schemes preserve the corresponding discrete dissipation law. Meanwhile, under some smoothness assumption, it is proved, by a higher order expansion technique, that both schemes are  second-order convergent in space and first-order convergent in time. Each scheme  yields a good approximation for the solution and the free boundary. No oscillation is observed for the numerical solution around the free boundary. Furthermore, the waiting time problem could be naturally treated, which has been a well-known difficult issue for all the existence methods. Due to its linear nature, the second scheme is more efficient.

 \end{abstract}
}

\section{Introduction and Background }
\label{sec:1}
The porous medium equation (PME) can be found  in many
  physical and biological phenomena, such as  the flow of an isentropic gas through a porous medium \cite{L. S. Leibenzon(1930)},  the viscous gravity currents \cite{J. Gratton(1998)}, nonlinear heat transfer and image processing; e.g., see \cite{J. L. Vazquez(2007)}. 
 The aim of  this paper is to provide numerical methods for the PME
 $$ \partial_t f=\Delta_{x} (f^{m}), \ x\in\Omega\subset\mathbb{R}^d,\  m>1,$$
where  $\ f:=f(x,t)$ is a non-negative scalar function of space $x\in\mathbb{R}^d$ and the time  $t\in\mathbb{R}$, the space dimension is given by $d\geq1$, and $m$ is a
constant larger than 1. 

The PME is a nonlinear degenerate parabolic equation since the diffusivity $D(f)=m f^{m-1}=0$ at points where $f = 0$. 
In turn, the PME has a special feature: the finite speed of  propagation, called \emph{finite} \emph{propagation} \cite{J. L. Vazquez(2007)}. If the initial data has a compact support, the solution of Cauchy problem of the PME will  have a compact support at any given time $t>0$. In comparison with the heat equation, which can smooth out the initial data, the solution of the PME becomes non-smooth even if the initial data is smooth with  compact support. If an initial data is zero in some open domain in $\Omega$, it causes the appearance of the free boundary (in some cases, called \emph{interface}) that separates the regions where the solution is positive from the regions where the value is zero in the domain.  Moreover, for certain initial data, the solution of the PME can exhibit a waiting time phenomenon where the free boundary remains stationary until a finite positive  time (called \emph{waiting time}). After that time instant, the interface begins to move with a finite speed.

Many  theoretical analyses have been available in the existing literature, including the earlier works by Ole\v{\i}nik et al. \cite{O.A. Oleinik(1958)}, Kala\v{s}nikov \cite{A.S. Kalasnikov(1967)}, Aronson \cite{D.G. Aronson(1969)}, the recent work by Shmarev \cite{S. Shmarev(2003),S. Shmarev(2005)} and the monograph by V\'azquez \cite{J. L. Vazquez(2007)}, etc. Among them, a fundamental example of solution is the Barenblatt solution \cite{G.I. Barenblatt (1952),R.E. Pattle (1959),Ya.B. Zeldovich (1950)}, which  has the explicit formula and 
   a compact
support at any  time  $t>0$ with   the interface.

Various numerical methods have been studied for the PME.  Graveleau \& Jamet \cite{J. L. Graveleau(1971)} and DiBenedetto \& Hoff \cite{E. DiBenedetto(1984)}  solved the pressure PME equation, using the finite difference approach  and tracking algorithm (containing a numerical viscosity term), respectively.   Jin et al. \cite{S. Jin(1998)}  established the relaxation scheme which reformulates the PME as a linear hyperbolic system with stiff relaxation term.  However, many existing numerical solutions may contain oscillations near the free boundary, such as PCSFE method (Predictor-Correction Algorithm and Standard Finite element method)  \cite{Q. Zhang(2009)}. 
 In recent years,  a  local discontinuous Galerkin finite element method by Zhang \& Wu  \cite{Q. Zhang(2009)}  and Variational Particle Scheme (VPS) by Westdickenberg \& Wilkening \cite{M. Westdickenberg(2010)} have been used to solve the PME. These two methods can effectively eliminate non-physical oscillation in the computed solution near the free boundary, and lead to a high-order convergence rate within the smooth part of the solution support.  However, no relevant theoretical justification of  the convergence analysis is available for these works.  More recently,  Huang \& Ngo \cite{C. Ngo(2017)} studied an adaptive moving mesh finite element method to solve the PME with three types of metric tensor: uniform, arclength-based and Hessian-based adaptive meshes. The numerical results indicate that a first-order convergence for uniform and arclength-based adaptive meshes, and a  second-order convergence for the Hessian-based adaptive mesh, while minor oscillations are observed around the free boundary in the computed solutions. Again, no theoretical proof has been available for the convergence rate in these works.
%

 For the  waiting phenomenon,  Mimura et al.  \cite{M. Mimura(1984)},   Bertsch  \& Dal Passo \cite{M. Bertsch (1990)}  and  Tomoeda \&  Mimura \cite{K. Tomoeda(1983)} estimated the waiting time  by the interface, but the numerical interface actually  has a velocity. Nakaki \& Tomoeda \cite{T. Nakaki(2003)} transformed the PME into another problem whose solution will blow up at a finite  time,  which is just the waiting time of PME. But the solution cannot be obtained  after the waiting time.
%



   In this paper, we construct numerical methods for PME by an Energetic Variational Approach (EnVarA) to naturally keep the physical laws, such as the conservation of mass, energy dissipation and force balance. Meanwhile, based on different dissipative energy laws, we can lead to different numerical schemes.
We start from the  energy dissipation law:
\begin{equation}\label{equ:EDLgen}
\frac{d}{dt}\int_{\Omega}\omega(f)dx=-\int_{\Omega}\eta(f)|{\bf u}|^2dx,
\end{equation}
where  $\omega(f)$ is the free energy density, $\eta(f)$ is a functional of $f$ determined by $\omega(f)$ and ${\bf u}$ is the velocity.  The quantity $\omega(f)$ and $\eta(f)$ can be taken as follows:
\begin{itemize}
\item Case 0. \ \ ${\displaystyle  \omega(f)=\frac{1}{m-1}f^m, \  \mbox{\ and\ } \eta(f)=f. }$

 \item Case 1. \ \ ${\displaystyle \omega(f)=f\ln f, \  \mbox{\ and\ } \eta(f)=\frac{f}{mf^{m-1}}. }$

   \item Case 2. \ \ ${\displaystyle  \omega(f)=\frac{1}{2f}, \ \mbox{\ and\ } \eta(f)=\frac{1}{mf^{m}}.} $
\end{itemize}
   Based on  these  energy dissipation laws, different  numerical schemes of the trajectory equation can be derived. The numerical scheme  based on the energy  law in Case 0  has been studied by Westdickenberg \& Wilkening  \cite{M. Westdickenberg(2010)},  called as a Variational Particle Scheme (VPS).

      We focus on the numerical methods based on the energy laws in the next two cases. Note that,  when $f$ vanishes, the energy in first case is regular while the energy in next two cases  is singular. Taking the advantage of the singularity, we can prove that the numerical schemes based on last two energy forms have some good properties which are not possessed by the VPS scheme \cite{M. Westdickenberg(2010)} from the first one, such as conservation of positivity,  unique solvability on an admissible convex set, convergence of the corresponding Newton's iteration.

 Theoretically, the discrete energy dissipation law is proved to be valid and by a higher order expansion technique  \cite{W. E(1995),C. Wang(2000)}, an optimal error estimates are derived under the assumption of smooth solutions.  Numerically, for Cases 1 and  2, no numerical oscillation is observed near the free boundary in the extensive   experiments, and the finite propagation speed of the free boundary can be effectively computed. A predictable criterion for computing waiting time is  proposed  and the numerical  convergence   to the exact waiting time is reported,  which is the first such result for  PME. In the practical computations, the numerical scheme of the trajectory equation in Case 2 is   linear  and hence more efficient.

   This paper is organized as follows. The EnVarA and the trajectory equation of the PME are outlined in Sec. \ref{sec:2}. The numerical scheme is described in  Sec. \ref{sec:3}. Subsequently, the proof of unique solvability, energy stability and optimal rate convergence analysis is provided in Sec. \ref{sec:4}. Finally, the numerical results are presented in Sec. \ref{sec:numerical results}, including examples with positive initial state,  Barenblatt Solution,  a waiting time phenomenon, an initial data with two columns, etc.

\section{Trajectory Equation of the PME}
\label{sec:2}



  In this section, we  derive the trajectory equation of the following initial-boundary problem of PME:
\begin{eqnarray}
&&
  \partial_t f=\Delta_{x} (f^{m}), \ x\in\Omega\subset\mathbb{R}^d,\  m>1, \ t>0,
 \label{eqPMEori}\\
&& f(x,0)=f_{0}(x)\geq 0, \ x\in\Omega, \label{eqPMEini} \\
&&\nabla_{x}f\cdot {\bf n}=0, \ x\in\partial\Omega,\ t>0, 
 \label{eqPMEboun}
 \end{eqnarray}
where $f$ is a non-negative function, $\Omega$ is a bounded domain and ${\bf n}$ is the external normal direction.
\subsection{The energetic variational approach}
An  Energetic Variational Approach (EnVarA) leads to the trajectory equation (also called \emph{constitution relation}) based on a balance between the maximal dissipation principle (MDP) and the least action principle (LAP). The approach was originated from Onsager's pioneering work \cite{L. Onsager(1931), L. Onsager1(1931)} and   improved by J.W. Strutt (Lord Rayleigh) \cite{J. W. Strutt(1873)}.  In recent years, it  has been applied to build up a mathematical model for a complex  physical  system, for example  Liu \&  Wu \cite{C. Liu(2003)},   Hyon et al. \cite{Y. Hyon(2010)}   Du et al. \cite{Q. Du(2009)},  Eisenberg et al. \cite{B. Eisenberg(2010)} and  Koba et al. \cite{Hajime Koba(2017)}. Its application to the Wright-Fisher model has been studied in  \cite{C.H. Duan(2017)}.The detailed structures of EnVarA can be found in \cite{C.H. Duan(2017),Y. Hyon(2010),C. Liu(2003),C. Liu(2017)}.

\noindent{\textbf{(A) Mass conservation.}}

In the Eulerian coordinate,  the mass conservation  law  is
 \begin{equation}\label{equ:conservationE}  \partial_t f+ \nabla\cdot(f\textbf{u})=0,\end{equation}
where $f$ is the density and $\textbf{u}$ is the velocity.

 In  the Lagrangian coordinate,  its solution can be expressed by:
 \begin{equation}\label{equ:conservationL}
  f(x(X,t),t)=\frac{f_{0}(X)}{\det \frac{\partial x(X,t)}{\partial X}},
\end{equation}
where $f_{0}(X)$ is the positive initial data  and  $\det \frac{\partial x(X,t)}{\partial X}$ is the  determinant of \emph{deformation\ gradient}.

\noindent{\bf(B) Energy Dissipation Law (EDL) Step.}

The basic energy dissipation law of PME we are going to consider is
\begin{equation}\label{eqEDL}
\frac{d}{dt}\int_{\Omega}\omega(f)dx=-\int_{\Omega}\eta(f)|{\bf u}|^2dx,
\end{equation}
where the total energy
  $E^{total}:=\int_{\Omega}\omega(f)dx$ with the free energy density $\omega(f)$,
and $\Delta:=\int_{\Omega}\eta(f)|{\bf u}|^2dx$ is the dissipation term with the  velocity ${\bf u}$.

  \noindent{\bf(C) Least Action  Principle (LAP) Step.}

 LAP states that   the trajectory of particles $X$ from the position $x(X,0)$ at time $t=0$ to
             $x(X,T^{*})$ at a given time $T^{*}$ in  Hamiltonian system are those which minimize the action functional defined by \ \
           $$\mathcal{A}(x):=-\int^{T^{*}}_{0}\mathcal{F}dt = -\int^{T^{*}}_{0}\int_{\Omega}\omega\left(\frac{f_0(X)}{\det\frac{\partial x}{\partial_X}}\right)\det\frac{\partial x}{\partial_X}dXdt,$$
           where  $\mathcal{F}$ is the Helmholtz free energy. 

          Taking the variational of $\mathcal{A}(x)$ with respect to $x$, we have the conservation force in Eulerian coordinate, i.e.,
            $$F_{con}:=\frac{\delta \mathcal{A}}{\delta x}=-\nabla(f\omega'(f)-\omega)=-f\nabla\omega'(f),$$
            where $\delta$ refers to the variational of the respective quantity.

  \noindent {\bf (D) Maximum Dissipation Principle (MDP) Step.}

      MDP i.e., Onsager's Principle, can be done by 
     taking the variational of $\frac{1}{2} \Delta$ with respect to the velocity $\textbf{u}$. In turn, we can obtain the dissipation force, i.e.,
       $$F_{dis}:=\frac{\delta\frac{1}{2}\Delta}{\delta\textbf{u}}=\eta(f){\bf u}.$$
       The factor $\frac{1}{2}$ is needed since that  the energy dissipation $\Delta$  is always a quadratic function of certain rates such as the velocity within the linear response theory \cite{J. W. Strutt(1873)}.

   \noindent   {\bf (E) Force  Balance  Law Step.}

      Based on the Newton's force balance  law:
       $$F_{con}=F_{dis},$$
       we have the constitution relation:
       $$f\nabla\omega'(f)=-\eta(f){\bf u}.$$
     which is just
     \begin{equation}\label{eqforce}
     \frac{f^2\omega''(f)\nabla f}{\eta(f)}=-f\bf u.
     \end{equation}
Comparing  PME \eqref{eqPMEori} with   \eqref{equ:conservationE}, we choose $-f\textbf{u}=\nabla(f^m)$, then

   $$\frac{f^2\omega''(f)}{\eta(f)}=mf^{m-1}.$$
    That means if  the free energy $\omega(f)$ is given, then $\eta(f)$ will be determined. Theoretically, there are infinite kinds of energy dissipation laws of PME. We consider three of them:
 \begin{itemize}
 \item
 Case 0.  if   $\omega(f)=\frac{1}{m-1}f^m$,  then $\eta(f)=f$ and the constitution  relation becomes
$$\nabla_x f^m=-f{\bf u}.$$
Let $P:=\frac{m}{m-1}f^{m-1}$ be  the pressure. The relation is the  Darcy's Law \cite{J. L. Vazquez(2007)}, i.e., $\textbf{u}=\nabla P$.

\item
 Case 1.  if $\omega(f)=f\ln f$, then $\eta(f)=\frac{f}{mf^{m-1}}$ and the constitution  relation in another form becomes
$$\nabla_x f= -\frac{f}{mf^{m-1}}\textbf{u}.$$

\item
 Case 2.  if $\omega (f)=\frac{1}{2f}$, then $\eta(f)=\frac{1}{mf^{m}}$ and the  constitution  relation in the third form is
$$\nabla_x \Big(\frac{1}{f}\Big)=\frac{\bf {u}}{mf^m}.$$
The free energy density $\frac{1}{2f}$ is a kind of  elastic energy \cite{T. Huang(2016)} and can lead to  a   linear numerical scheme for the trajectory equation.
\end{itemize}

\subsection{ Trajectory Equation  of PME in 1-Dim}
Combining with   \eqref{equ:conservationL}, we can  write the constitution relation  in the Lagrangian coordinate system, called as the trajectory equation. In this paper, we consider one dimensional problems. Replacing ${\bf u}$ with $x_t(X,t)$, we have the trajectory equation as
 \begin{itemize}
\item
Case 0.  
\begin{equation}\label{eqtra1}f_0(X)\partial_t x=-\partial_X\left[\Big(\frac{f_0(X)}{\partial_X x}\Big)^m\right],\  \ X\in\Omega,\end{equation}
 and   the corresponding energy law in  Lagrangian coordinate is
  $$\frac{d}{dt}\int_{\Omega}\frac{1}{m-1}\Big(\frac{f_0(X)}{\partial_X x}\Big)^m\frac{\partial x}{\partial X} dX=-\int_{\Omega}f_0(X)|x_t|^{2}dX.$$
  \item
   Case 1.
 \begin{equation}\label{eqtra2}\frac{f_{0}(X)}{m\big(\frac{f_0(X)}{\partial_X x}\big)^{m-1}}\partial_{t} {x}=-\partial_{X}\left(\frac{f_0(X)}{\partial_X x}\right), \ \ X\in\Omega,\end{equation}
 and   the corresponding energy law in  Lagrangian coordinate is
 \begin{equation} \label{tot-energ-e2} \frac{d}{dt}\int_{\Omega}f_0(X)\ln \Big(\frac{f_0(X)}{\partial_X x}\Big)dX=-\int_{\Omega}\frac{f_0(X)}{m\big(\frac{f_0(X)}{\partial_X x}\big)^{m-1}}|\partial_{t} {x}|^{2}dX,\end{equation}
 \item
  Case 2. 
\begin{equation}\label{eqtra3}\frac{(\partial_X x)^{m+1}}{mf_0(X)^m}\partial_t x=\partial_X\Big(\frac{\partial_X x}{f_0(X)}\Big),\  \ X\in\Omega,\end{equation}
and   the corresponding energy law in  Lagrangian coordinate is
   \begin{equation}\label{tot-energ-e3}\frac{d}{dt}\frac{1}{2}\int_{\Omega}\frac{1}{f_0(X)} |\partial_X x|^2dX=-\int_{\Omega}\frac{\partial_X x}{m\Big(\frac{f_0(X)}{\partial_X x}\Big)^{m}}|x_t|^{2}dX.\end{equation}
\end{itemize}
 Equations \eqref{eqtra1},  \eqref{eqtra2} and \eqref{eqtra3} are the same thing, wrote  in different forms with differnet energy laws. Solving them with proper initial and boundary conditions, we get the  trajectory $x(X,t)$, which contains all the physics involved in the model.  Substituting  $x(X,t)$   into \eqref{equ:conservationL}, we obtain the   solution $f(x,t)$ to \eqref{eqPMEori}-\eqref{eqPMEboun}.

The initial and boundary conditions for equations \eqref{eqtra1}, \eqref{eqtra2} or \eqref{eqtra3} should be
\begin{equation}\label{eqtraboun} x|_{\partial\Omega}= X|_{\partial\Omega},\ t>0. \end{equation}
\begin{equation}\label{eqtraini}  x(X,0)=X, \ X\in\Omega. \end{equation}



\section{Numerical Method of Trajectory Equation}
\label{sec:3}
In this section, we propose some  semi-implicit numerical schemes  for  the trajectory equations.

\subsection{Semi-discrete schemes in time}

 Let $\tau=\frac{T}{N}$, where $N\in\mathbb{N}^+$ and  $T$  is the final time. The grid point $t_n=n\tau$, $n=0,\cdots,N$.
For  the temporal discretization of the trajectory equation, we have that Given $x^n$, find $x^{n+1}$ such that
   \begin{itemize}
   \item
   Case 0.  
      \begin{equation}\label{eqtradnum1}
     f_{0}(X)\frac{x^{n+1}-x^{n}}{\tau}=- \partial_X\left[\Big(\frac{f_0(X)}{\partial_X x^{n+1}}\Big)^m\right],\  n=0,\cdots,N-1.
    \end{equation}

   \item
     Case 1.  
    \begin{equation}\label{eqtradnum2}
     \frac{f_{0}(X)}{m\Big(\frac{f_0(X)}{\partial_X x^{n}}\Big)^{m-1}}\frac{x^{n+1}-x^{n}}{\tau}=- \partial_X\Big(\frac{f_0(X)}{\partial_X x^{n+1}}\Big),\  n=0,\cdots,N-1.
    \end{equation}

    \item
   Case 2. 
    \begin{equation}\label{eqtradnum3}
     \frac{(\partial_X x^{n})^{m+1}}{mf_0(X)^{m}}\frac{x^{n+1}-x^{n}}{\tau}= \partial_X\Big(\frac{\partial_X x^{n+1}}{f_0(X)}\Big),\  n=0,\cdots,N-1.
    \end{equation}
\end{itemize}
In summary,  the trajectory equation can be written in gradient flow as
\begin{equation}\label{equ:traGF}\gamma(x)x_t=-\frac{\delta\mathcal{W}}{\delta x},\end{equation}
where $\gamma(x,t)$ is a positive function depending on space $x$ and time $t$ and $\mathcal{W}$ is a functional of $x$. Then the  discrete scheme in time is
\begin{equation}\label{equ:traGFTime}
\gamma(x^n)\frac{x^{n+1}-x^{n}}{\tau}=-\frac{\delta\mathcal{W}(x^{n+1})}{\delta x^{n+1}}, \ n=0,\cdots, N-1.
\end{equation}

  Assume the exact solution $x^n$ is smooth at time $t^n$ to make $\frac{\partial x}{\partial X}$ well-defined, $n=0,\cdots,N$, then the solution $x^{n+1}$ to the numerical scheme  to \eqref{equ:traGFTime} is   the minimizer of the following cost functional:
  \begin{equation}\label{equ:traGFOptimal}\min\limits_{x^{n+1}\in\Omega}\left\{\int_{\Omega}\gamma(x^n)\frac{|x^{n+1}-x^{n}|^2}{2\tau}+\mathcal{W}(x^{n+1})dX\right\},\end{equation}
  where
  \begin{itemize}
  \item
 Case 0
   $$\gamma(x^n)=1,\ \mathcal{W}(x^{n+1})=\frac{1}{m-1}\frac{f_0(X)^{m}}{(\partial_X x^{n+1})^{m-1}},$$ where $f_0(X)$ is the initial function.
  It  means that the trajectory equation can be regarded as an energy gradient flow, which has been studied by Westdickenberg and Wilkening  \cite{M. Westdickenberg(2010)}. 
   In this paper, we focus on the following two cases:

 \item
  Case 1 
  $$\gamma(x^n)=\frac{f_{0}(X)}{m\big(\frac{f_0(X)}{\partial_X x^{n}}\big)^{m-1}},\
\mathcal{W}(x^{n+1}) =f_0(X)\ln\Big(\frac{f_0(X)}{\partial_X x^{n+1}}\Big).$$

\item
  Case 2 
  $$\gamma(x^n)=\frac{(\partial_X x^{n})^{m+1}}{mf_0(X)^{m}},\ \mathcal{W}(x^{n+1}) = -\frac{1}{2}\frac{1}{f_0(X)}|\partial_X x^{n+1}|^2.$$
\end{itemize}
\subsection{The fully discrete scheme with a positive initial state}
Let $X_0$ be the left point of $\Omega$ and $h=\frac{|\Omega|}{M}$  be the spatial step,  $M\in\mathbb{N}^{+}$. Denote by $X_{r}=X(r)=X_0+ r h$, where $r$ takes on integer and half integer values.  Let $\mathcal{E}_{M}$ and $\mathcal{C}_{M}$ be the spaces of functions whose domains are $\{X_{i}\ |\ i=0,...,M\}$ and $\{X_{i-\frac{1}{2}}\ |\ i=1,...,M\}$, respectively. In component form, these functions are identified via
$l_{i}=l(X_{i})$, $i=0,...,M$, for $l\in\mathcal{E}_{M}$, and $\phi_{i-\frac{1}{2}}=\phi(X_{i-\frac{1}{2}})$,  $i=1,...,M$, for $\phi\in\mathcal{C}_{M}$.

\indent{The} difference operator $D_{h}: \mathcal{E}_{M}\rightarrow\mathcal{C}_{M}$, $d_{h}: \mathcal{C}_{M}\rightarrow\mathcal{E}_{M}$, and $\widetilde{D}_h: \mathcal{E}_{M}\rightarrow\mathcal{E}_{M}$  can be defined  as:
\begin{align}\label{equ:dif1}
& (D_{h}l)_{i-\frac{1}{2}}= (l_{i}-l_{i-1})/h,\ i=1,...,M, \\
&  (d_{h}\phi)_{i}= (\phi_{i+\frac{1}{2}}-\phi_{i-\frac{1}{2}})/h,\ i=1,...,M-1,\\
&(\widetilde{D}_hl)_{i}=\left\{
\begin{array}{lcl}
(l_{i+1}-l_{i-1})/2h, &\mbox{$ i=1,...,M-1$},\\
 (4l_{i+1}-l_{i+2}-3l_{i})/2h,&  i=0,  \\
 (l_{i-2}-4l_{i-1}+3l_{i})/2h, &  i=M.\end{array}\right.
\end{align}

Let $\mathcal{Q}:=\{l \in\mathcal{E}_{M}\ |\ l_{i-1}<l_{i},\ 1\leq i\leq M;\ l_{0}=X_0,\ l_{M}=X_M\}$ be the admissible set, in which the particles are arranged in the order without twisting or exchanging. Its boundary set is $\partial\mathcal{Q}:=\{l \in\mathcal{E}_{M}\ |\ l_{i-1}\leq l_i,\ 1\leq i\leq M,\ and\  l_{i}=l_{i-1},\ for\ some\ 1\leq i\leq M;\ l_{0}=X_0,\ l_{M}=X_M\}$. Then $\bar{\mathcal{Q}}:=\mathcal{Q}\cup\partial\mathcal{Q}$ is a closed convex set.

 \indent The {\bf  fully discrete scheme} is formulated as follows. Given the positive initial state $f_0(X)\in\mathcal{E}_M$ and the particle position $x^{n} \in\mathcal{Q}$, find $x^{n+1}=(x^{n+1}_{0},...,x^{n+1}_{M})\in\mathcal{Q}$ such that


\begin{itemize}
\item
  Case 1.  
  \begin{eqnarray}
  &&\label{eqtranum2}
\frac{f_{0}(X_{i})}{m\big(\frac{f_0(X)}{\widetilde{D}_h x^{n}}\big)_{i}^{m-1}}\cdot\frac{x^{n+1}_{i}-x^{n}_{i}}{\tau}=-d_{h}\Big(\frac{f_{0}(X)}{D_{h} x^{n+1}}\Big)_{i}, \ 1\leq i \leq M-1,\\
&& \mbox{with\ \ } x_0^{n+1}=X_0  \mbox{\ and\ } x_M^{n+1}=X_M, \ n=0,\ldots,N-1. \nonumber
\end{eqnarray}

To solve the nonlinear equation \eqref{eqtranum2}, we use damped Newton's iteration\cite{Y. Nesterov(1994)}.  The key idea is to adjust the marching size to prevent the solution at next iteration to escape from the admissible set $\mathcal{Q}$.

\noindent{\bf Damped Newton's iteration}. \ \ Set $x^{n+1, 0}= x^n$. For $k=0,1,2,\cdots$, update $x^{n+1,k+1}=x^{n+1,k} + \omega(\lambda)\delta_x$ such that
\begin{eqnarray}\label{equ:numnum1}
&&
\frac{f_{0}(X_{i})}{m\big(\frac{f_0(X)}{\widetilde{D}_h x^{n}}\big)_{i}^{m-1}}\frac{ \delta_{x_{i}}}{\tau} - d_h \left(\frac{f_{0}(X)}{(D_{h}x^{n+1,k})^2} D_h \delta_{x}\right)_{i}\nonumber\\
&&= -
 \frac{f_{0}(X_{i})}{m\big(\frac{f_0(X)}{\widetilde{D}_h x^{n}}\big)_{i}^{m-1}}\frac{x^{n+1,k}_{i} -x^{n}_{i}}{\tau} -d_{h}\left(\frac{f_{0}(X)}{D_{h}x^{n+1,k}}\right)_{i},
\ \ 1\leq i\leq M-1,\\
&&
\mbox{with\ \ }\delta_{x_0}=\delta_{x_M}=0,  \nonumber
\end{eqnarray}
and
\begin{equation}\label{Newton}
\omega(\lambda)=\left\{
\begin{array}{lcl}
\frac{1}{\lambda}, &&{\lambda > \lambda'},\\
\frac{1-\lambda}{\lambda(3-\lambda)}, &&{\lambda'\ge\lambda\ge\lambda^{*}},\\
1, && {\lambda < \lambda^{*}},
\end{array}\right.
\end{equation}
where $\lambda^*=2-3^{\frac{1}{2}}$, $\lambda'\in[\lambda^*, 1)$ and
\begin{equation}\label{lambda}
\lambda^2:=\lambda^2(J,x^{n+1,k})=\frac{1}{a}[J'(x^{n+1,k})]^T[J''(x^{n+1,k})]^{-1}J'(x^{n+1,k}),
\end{equation}
where $a:=h\min\limits_{0< i < M}\{f_0(X_{i})\}$,  $J$ is the corresponding energy function  defined latter in \eqref{DisEnergy}, and $J', J''$ are the gradient vector and Hessian matrix.
 \item
  Case 2.
 \begin{eqnarray}
 &&\label{eqtranum3}
 \frac{(\widetilde{D}_h  x^n )^{m+1}_i}{mf_0(X_i)^{m}}\cdot\frac{x^{n+1}_{i}-x^{n}_{i}}{\tau}=d_{h}
\Big(\frac{D_h x^{n+1}}{f_0(X)}\Big)_{i}, \ 1\leq i \leq M-1,\\
 && \mbox{with\ } x_0^{n+1}=X_0 \mbox{\ and\ } x_M^{n+1}=X_M,  \ n=0,\cdots,N-1.\nonumber
 \end{eqnarray}
Note that \eqref{eqtranum3} is  a linear  scheme.
\end{itemize}

After solving \eqref{eqtranum2} in Case 1  (\eqref{eqtranum3} in Case 2),  we finally obtain the numerical solution $f_{i}^{n}:=f(x^{n},t^{n})$ from \eqref{equ:conservationL} by
\begin{equation}\label{numdist}
f_{i}^{n}= \frac{ f_{0}(X_{i}) } {\widetilde{D}_h x^{n}_{i}},\ 0\leq i\leq M.
\end{equation}

\subsection{The discrete scheme  for problems with free boundaries}
Next we consider the situation of the initial data with a compact support in $\Omega$.  Due to the degeneration of the PME, the left and right  interfaces appear and are defined respectively as:
$$\xi_{1}^t:=\inf\{x\in\Omega:f(x,t)>0,t\geq 0\},$$
$$\xi_{2}^t:=\sup\{x\in\Omega:f(x,t)>0,t\geq 0\}.$$
 Let $\Gamma^{t}:=[\xi_{1}^t, \xi_{2}^t]\subset\Omega$. For this kind of problems, all the trajectories start from the initial support $\Gamma^0\subsetneqq\Omega$. We shall solve a initial-boundary  value problem as: 
\begin{itemize}
\item
Case 1. 
 \begin{eqnarray}
 &&    \frac{f_{0}(X)}{m\big(\frac{f_0(X)}{\partial_X x}\big)^{m-1}}\partial_{t} {x}=-\partial_{X}\left(\frac{f_0(X)}{\partial_X x}\right),\ X\in\Gamma^0,\ t>0,  \label{eqtraC2}\\
 && (\partial_X x)^{m-1}\cdot \partial_t x=-\frac{m}{m-1}\frac{\partial_X[f_0(X)^{m-1}]}{\partial_X x},\ X\in\partial\Gamma^0, \  t>0, \label{eqbouC2}\\
     && x(X,0)=X, \ X\in\Gamma^0,  \label{eqiniC2}
    \end{eqnarray}
\item    Case 2.
       \begin{eqnarray}
 &&  \label{eqtraC3}\frac{(\partial_X x)^{m+1}}{mf_0(X)^m}\partial_t x=\partial_X\Big(\frac{\partial_X x}{f_0(X)}\Big), \ \ X\in\Gamma^0,\ t>0,  \\
 && (\partial_X x)^{m-1}\cdot \partial_t x=-\frac{m}{m-1}\frac{\partial_X[f_0(X)^{m-1}]}{\partial_X x},\ X\in\partial\Gamma^0, \  t>0, \label{eqbouC3}\\
     && x(X,0)=X, \ X\in\Gamma^0.  \label{eqiniC3}
    \end{eqnarray}
    \end{itemize}
\begin{rem}
Taking into account that $f_0(X) = 0$ at the boundary of its support,  the boundary equation \eqref{eqbouC2} and \eqref{eqbouC3} is just the equation \eqref{eqtraC2}.
\end{rem}
Let   $h:=(\xi_{2}^0-\xi_{1}^0)/M$ be  the spatial step.
Then we partition the interval $\Gamma^0$ into equal subinterval with $X_i=\xi_{1}^0+ih$, $0\leq i \leq M$.

The \textbf{fully discrete} scheme becomes: Given  the initial state $f_0(X)$ with a compact support $\Gamma^0$ and $\{x_i^n\}_{i=0}^M$, find $\{x_i^{n+1}\}_{i=0}^M$ such that


\begin{itemize}
\item     Case 1. 
   \begin{eqnarray} && \label{equtrafreenum2}
\frac{f_{0}(X_{i})}{m\big(\frac{f_0(X)}{\widetilde{D}_h x^{n}}\big)_{i}^{m-1}}\cdot\frac{x^{n+1}_{i}-x^{n}_{i}}{\tau}=-d_{h}\Big(\frac{f_{0}(X)}{D_{h} x^{n+1}}\Big)_{i}, \ 0< i < M, \\
&& \label{bouF}  (\bar{D}_h x^n_i)^{m-1}\cdot\frac{x^{n+1}_i-x^{n}_i}{\tau}=-\frac{m}{m-1}\cdot\frac{\bar{D}_h [f_0(X_i)^{m-1}]}
 {\bar{D}_h x^{n+1}_i}, \label{equtrafreenumboun2}\ \ i=0,M,
 \end{eqnarray}
where
\begin{equation}\label{deri1}
\bar{D}_h l_i:=\left\{
\begin{array}{lcl}
\frac{l_{i+1}-l_{i}}{h}, && i=0,\\

\frac{l_{i}-l_{i-1}}{h}, && i=M,
\end{array}\right.   \quad \forall l=(l_0,\cdots,l_M).\end{equation}
 \item   Case 2.
 \end{itemize}
 \begin{eqnarray}
 &&
\frac{(\widetilde{D}_h  x^n )^{m+1}_i}{mf_0(X_i)^{m}}\cdot\frac{x^{n+1}_{i}-x^{n}_{i}}{\tau}=
d_{h}\Big(\frac{D_h x^{n+1}}{f_0(X)}\Big)_{i}, \ 1\leq i \leq M-1,\label{equtrafreenum3}\\
&&
(\bar{D}_h x^n_i)^{m-1}\cdot\frac{x^{n+1}_i-x^{n}_i}{\tau}=-\frac{m}{m-1}\cdot
\frac{\bar{D}_h [f_0(X_i)^{m-1}]}
 {\bar{D}_h x^{n+1}_i}, \ i=0,M.\label{equtrafreenumboun3} 
\end{eqnarray}

Comparing with the schemes \eqref{eqtranum2} and \eqref{eqtranum3}, we have two more nonlinear equations at the boundary.  The damped Newton's iteration shall be applied to solve the whole system.

\begin{rem} Note that the equation \eqref{equtrafreenum3} is linear but the boundary equation \eqref{equtrafreenumboun3} is chosen to be nonlinear, the same as \eqref{equtrafreenumboun2}. If we choose a linear boundary equation  as
$$ (\bar{D}_h x^n_i)^{m+1}\cdot\frac{x^{n+1}_i-x^{n}_i}{\tau}=-\frac{m}{m-1}\cdot
 \bar{D}_h [f_0(X_i)^{m-1}] \cdot \bar{D}_h x^{n+1}_i, \ i=0,M,$$
then the matrix of the whole linear system would not be a M-matrix and the conservation of positivity would be destroyed.
\end{rem}

When the right side of equation \eqref{equtrafreenumboun2} is zero,  the waiting  phenomenon occurs. During the waiting time, the boundary condition in \eqref{eqbouC2} or \eqref{eqbouC3} should be replaced by $x_t = 0,  X\in\partial\Gamma^0$ and the boundary condition in \eqref{equtrafreenumboun2} or \eqref{equtrafreenumboun3} should be replaced by $x^{n+1}_0 = \xi_1^0,\  x^{n+1}_M = \xi_2^0$. The key problem is how to predict when the waiting stops.  For the details to treat this kind of problem, see the algorithm in Section \ref{sec:numerical results}, Example 3.

\section{Analysis of the Numerical Schemes}
\label{sec:4}
In this section, we perform detailed analyses for the numerical schemes \eqref{eqtranum2} and  \eqref{eqtranum3} , including the unique solvability in admissible  set, the optimal rate convergence analysis, the convergence of Newton's iteration and  the dissipation analysis of the total energy. 

A few more notations have to be introduced.  Let $l$, $g\in\mathcal{E}_M$ and $\phi$, $\varphi\in\mathcal{C}_M$. We define the \emph{inner product} on space $\mathcal{E}_M$ and $\mathcal{C}_M$ respectively as:
\begin{eqnarray}
  &&
 \left\langle l , g \right\rangle  := h \left( \frac12 l_0 g_0 + \sum_{i=1}^{M-1}   l_{i} g_{i}
  + \frac12 l_M g_M \right),\\
  \label{FD-inner product-1}
&&
  \left\langle \phi , \varphi \right\rangle_e
  := h \sum_{i=0}^{M-1}   \phi_{i+\frac 12} \varphi_{i+\frac 12} .
  \label{FD-inner product-2}
\end{eqnarray}

The following summation by parts formula is available:
\begin{equation}
   \left\langle l , d_h \phi \right\rangle
  = -\left\langle D_h l , \phi \right\rangle_e,
  \mbox{\ with $l_0 = l_M =0$, $\phi\in\mathcal{C}_{M}$, $l\in\mathcal{E}_M$}.
  \label{FD-inner product-3-1}
\end{equation}

The inverse inequality is available:
\begin{equation}\label{equ:inverse}
\|l\|_{\infty}\leq C_m \frac{\|l\|_2}{h^{1/2}}, \ \ \forall l\in\mathcal{E}_M,
\end{equation}
where $$\|l\|_{\infty}:=\max\limits_{0\leq i\leq M} \{l_i\}\mbox{\ \ and\ \ }\|l\|_2^2:=\left\langle l,l \right\rangle.$$

First we prove that there exists a unique solution in admissible  set $\mathcal{Q}$.

\begin{thm}
\label{lem:unique}
Suppose $f_0(X)\in\mathcal{E}_{M}$ is the initial state  with a positive lower bound for $X\in \mathcal{Q}$. The numerical scheme  \eqref{eqtranum2} is uniquely solvable in  $\mathcal{Q}$, and the solution $x^{n+1}$ to the linear  scheme  \eqref{eqtranum3} also belongs to $\mathcal{Q}$, for $n=1,\cdots,N-1$. 
\end{thm}
\noindent\textbf{Proof:} To prove the existence and uniqueness of solution in $\mathcal{Q}$ to the scheme \eqref{eqtranum2},  we first consider the following optimization problem:
\begin{equation}\label{DisEnergy}
\min\limits_{y\in\bar{\mathcal{Q}}} J(y):=\frac{1}{2\tau}\Big\langle \frac{f_{0}(X)}{m\big(\frac{f_0(X)}{\widetilde{D}_h x^n}\big)^{m-1}}
(y-x^{n}),(y-x^{n})\Big\rangle+
\Big\langle f_{0}(X),\ln\Big(\frac{ f_0(X)}{D_h y}\Big)\Big\rangle_e,
\end{equation}
where  $x^{n}\in\mathcal{Q}$ is the   position of particles at time $t^{n}$, $n=0,\cdots,N-1$. Since $J(y)$ is a convex function on the closed convex set $\bar{\mathcal{Q}}$,  there exists a unique minimizer  $x\in\bar{\mathcal{Q}}$. Moreover,  we must have  $x\in\mathcal{Q}$, since for $\forall$ $y\in\partial\mathcal{Q}$,  there exists some $i>0$ such that $(D_h y)_{i-1/2} = (y_i-y_{i-1})/h=0$, then $J(y)=+\infty$.

Next we  want to prove  that $x\in\mathcal{Q}$ is the minimizer of $J(y)$  if and only if it is a solution to   scheme \eqref{eqtranum2}. Then we can claim that  the fully discrete scheme \eqref{eqtranum2}  has a unique solution.

In fact, if $x\in\mathcal{Q}$ is the minimizer of $J(y)$, then for $\forall y\in\bar{\mathcal{Q}}$, there exists a sufficiently small $\varrho_0>0$, such that for any $\varrho\in(-\varrho_0, \varrho_0)$, $x+\varrho (y-x) \in \mathcal{Q}$ since $\mathcal{Q}$ is a open convex set. Then $j(\varrho) := J(x+\varrho (y-x))$ achieves its minimal at $\varrho =0$.  So we have $j'(0)=0$ and  using summation by parts, we obtain

\begin{equation*}
    \frac{1}{\tau}\Big\langle \frac{f_{0}(X)}{m\big(\frac{f_0(X)}{\widetilde{D}_h x^n}\big)^{m-1}}(x-x^{n}),y-x\Big\rangle
    +\Big\langle d_h\Big(\frac{f_0(X)}{D_h x^{n+1}}\Big), y-x\Big\rangle=0,
\end{equation*}
for any $y\in\bar{\mathcal{Q}}$. This implies that  $x\in\mathcal{Q}$ satisfies \eqref{eqtranum2}.

\indent{Conversely} let $x\in\mathcal{Q}$ be the solution to scheme \eqref{eqtranum2}. We need to prove that $x$ is the minimizer of $J(y)$ on $ \bar{\mathcal{Q}}$.

 For any  $ y\in\partial\mathcal{Q}$, we always have  $J(y)\geq J(x)$ due to $J(y)=+\infty$.
 Then for  any $y\in\mathcal{Q}$, taking the inner product  of (\ref{eqtranum2}) with $y-x$ and using summation by parts, we have
\begin{equation}\label{equ:vp1}
\frac{1}{\tau}\Big\langle \frac{f_{0}(X)}{m\big(\frac{f_0(X)}{\widetilde{D}_h x^n}\big)^{m-1}}(x-x^{n}),y-x\Big\rangle
    -\Big\langle \frac{f_0(X)}{D_h x}, D_h(y-x)\Big\rangle_e=0.
    \end{equation}
After the direct calculation,  we  get  for  any $y\in\mathcal{Q}$ such that
  \begin{align}
  J(y)&=J(x+(y-x)) =J(x)+\frac{1}{2\tau}\Big\langle \frac{f_{0}(X)}{m\big(\frac{f_0(X)}{\widetilde{D}_h x^n}\big)^{m-1}}(y-x),y-x\Big\rangle
     \nonumber\\
     &\ \ \ \ +\frac{1}{\tau}\Big\langle \frac{f_{0}(X)}{m\big(\frac{f_0(X)}{\widetilde{D}_h x^n}\big)^{m-1}}(x-x^{n}),y-x\Big\rangle+\Big\langle f_{0}(X), \ln\Big(\frac{D_h x}{D_h y}\Big)\Big\rangle_e \nonumber  \\
  &\geq J(x), \label{equ:vp3}
  \end{align}
 where the last inequality is obtained from  (\ref{equ:vp1})   and the fact: $\ln \frac{1}{z} \geq -(z-1)$ ,  $\forall z\in\mathbb{R}^+$, which leads to
 $$\Big\langle f_{0}(X), \ln\Big(\frac{D_h x}{D_h y}\Big)\Big\rangle_e\geq -\Big\langle f_{0}(X), \frac{D_h(y-x)}{D_h x}\Big\rangle_e.$$

Then we prove that  the solution to the numerical scheme \eqref{eqtranum3} $x^{n+1}\in\mathcal{Q}$ if given $x^n\in\mathcal{Q}$, $n=0,\cdots,N-1$.
Without loss of generality, let $\bar\Omega=[0,1]$.  Due to  the boundary condition \eqref{eqtraboun}, we have $$x_0^{n+1}=0,\ x_M^{n+1}=1.$$
Based on the discrete extremum principle, we obtain that
\begin{equation}\label{equ:PMELthm:unique1}
 0< x^{n+1}_i<1, \ i=1,\cdots,M-1,\ n=0,\cdots,N-1.
\end{equation}
Suppose $x^{n+1}\notin\mathcal{Q}$,  i.e., $\exists\  k_1,\ k_2\in\mathbb{N}^+$ such that $k_1<k_2$ and
\begin{equation}\label{equ:PMELthm:unique2}
0<x_{k_1-1}^{n+1}< x_{k_1}^{n+1}\geq x_{k_1+1}^{n+1} \geq \cdots \geq x_{k_2-1}^{n+1}\geq x_{k_2}^{n+1}<x_{k_2+1}^{n+1}<1.
\end{equation}
Checking the equation \eqref{eqtranum3} at $i=k_1$ and $i=k_2$ respectively,   we have
$$x_{k_1}^{n+1}<x_{k_1}^{n}<x_{k_2}^{n}<x_{k_2}^{n+1},$$
which contradicts with \eqref{equ:PMELthm:unique2}. Due to the initial state $X\in\mathcal{Q}$, then $x^{n}\in\mathcal{Q}$, $n=0,\cdots,N$.
The proof is finished. $\hfill\Box$

Next we  prove that the numerical scheme \eqref{eqtranum2}  and  \eqref{eqtranum3} satisfy  the corresponding  discrete  energy dissipation laws. 

\begin{thm}
Suppose the initial state $f_0(X)\in\mathcal{E}_{M}$ is   positive and bounded for $X\in\mathcal{Q}$.
\begin{itemize}
\item Case 1. Let $x^{n}=(x_{0}^{n},...,x^{n}_{M}) \in \mathcal{Q}$, $n=0,1,\cdots, N-1$, be the solution to scheme \eqref{eqtranum2} at time $t^{n}$. Then   the discrete  energy dissipation law holds, i.e.,
\begin{equation}\label{PMELequ:disenergy2}\frac{E^{(1)}_{N}(x^{n+1})-E^{(1)}_{N}(x^n)}{\tau } \leq - \Big\langle\frac{f_{0}(X)}{m\big(\frac{f_0(X)}{\widetilde{D}_h x^n}\big)^{m-1}}\frac{x^{n+1}-x^{n}}
{\tau}, \frac{x^{n+1}-x^{n}}
{\tau}\Big\rangle,
\end{equation}
where \begin{equation}  E^{(1)}_{N}(x):=\Big\langle f_{0}(X),\ln\Big(\frac{f_{0}(X)}{D_{h}x}\Big)\Big\rangle_e, \ \ \mbox{with\ }
  \frac{\delta E_{N}^{(1)}(x)}{\delta x}=d_{h}\Big(\frac{f_{0}(X)}{D_{h}x}\Big).\label{variation2}\end{equation}
\item Case 2. Let $x^{n}=(x_{0}^{n},...,x^{n}_{M}) \in \mathcal{Q}$, $n=0,1,\cdots, N-1$, be the solution to scheme \eqref{eqtranum3} at time $t^{n}$. Then   the following discrete energy dissipation law holds, i.e.,
\begin{equation}\label{PMELequ:disenergy3}\frac{E^{(2)}_{N}(x^{n+1})-E^{(2)}_{N}(x^n)}{\tau } \leq - \Big\langle\frac{(\widetilde{D}_h  x^n )^{m+1}}{mf_0(X)^{m}}\cdot\frac{x^{n+1}-x^{n}}{\tau}, \frac{x^{n+1}-x^{n}}
{\tau}\Big\rangle,  \end{equation}
where \begin{equation} E^{(2)}_{N}(x):=\frac{1}{2}\Big\langle \frac{D_{h}x}{f_{0}(X)},D_{h}x\Big\rangle_e,\ \ \mbox{with\ }
\frac{\delta E_{N}^{(2)}(x)}{\delta x}=-d_{h}\Big(\frac{D_{h}x}{f_{0}(X)}\Big).  \label{variation3} \end{equation}
\end{itemize}
\end{thm}
Note that \eqref{PMELequ:disenergy2} and \eqref{PMELequ:disenergy3} are  the  discrete counterpart of energy laws \eqref{tot-energ-e2} and  \eqref{tot-energ-e3}.

\noindent{\textbf{Proof.}}
In Case 1, thanks to the convexity of $E_{N}^{(1)}(x)$, we have
\begin{align}
\frac{E_{N}^{(1)}(x^{n})-E_{N}^{(1)}(x^{n+1})}{\tau}&\geq
  \Big\langle\frac{\delta E_{N}^{(1)}(x^{n+1})}{\delta x},\frac{x^{n}-x^{n+1}}{\tau}\Big\rangle =\Big\langle d_{h}\big(\frac{f_{0}(X)}{D_{h}x^{n+1}}\big), \frac{x^{n}-x^{n+1}}{\tau}\Big\rangle\notag\\
  &=\Big\langle\frac{f_{0}(X)}{m\big(\frac{f_0(X)}{\widetilde{D}_h x^n}\big)^{m-1}}\frac{x^{n}-x^{n+1}}{\tau},\frac{x^{n}-x^{n+1}}{\tau}\Big\rangle.\notag
\end{align}
That means \eqref{PMELequ:disenergy2} holds. Due to the convexity of $E_{N}^{(2)}(x)$, we can also  prove that the numerical scheme   \eqref{eqtranum3}  satisfies the discrete energy dissipation law \eqref{PMELequ:disenergy3} in the similar way.
  $\hfill\Box$

Next we provide the optimal rate convergence analysis for the schemes \eqref{eqtranum2}  and  \eqref{eqtranum3}.
\begin{thm}
\label{convergence}
Assume that  the initial function $f_0(X)$ is positive and bounded, i.e., $0<b_f \leq f_0(X) \leq B_f$.
 Denote $x_e\in\Omega$ as the exact solution to the original trajectory equation \eqref{eqtra2} or \eqref{eqtra3} with enough regularity and $x_{h}\in\mathcal{Q}$ as the numerical solution to  the numerical scheme \eqref{eqtranum2} in Case 1 or \eqref{eqtranum3} in Case 2. The numerical error function is defined at a point-wise level:
\begin{equation}\label{error-function-1}
 e_i^n = x_{e_i}^n - x_{h_i}^n ,
\end{equation}
where $x_{e_i}^n,\ x_{h_i}^n\in\mathcal{Q}$, $0 \le i \le N$, $n=0,\cdots,M$. 
Then
\begin{itemize}
\item{
$e^{n}=(e^{n}_0,\cdots,e^{n}_M)$ satisfies}
$$\quad  \| e^{n}\|_2:=\langle e^{n},e^{n}\rangle\le C (\tau + h^2).$$
\item{$\widetilde{D}_h e^{n}=(\widetilde{D}_h e^{n}_0,\cdots,\widetilde{D}_h e^{n}_M)$ satisfies}
$$\quad  \| \widetilde{D}_h e^{n} \|_{2} \leq C(\tau+h^2).$$
Moreover, the error between the numerical solution $f_h^{n}$ and the exact solution $f_e^{n}$ of the  problem \eqref{eqPMEori}-\eqref{eqPMEboun} can be estimated by: 
$$\quad  \|f_h^{n}-f_e^{n}\|_2 \leq C(\tau+h^2),$$
where $C$ is a positive constant, $h$ is the spatial step, $\tau$ is the time step and $n=0,\cdots,N$.
\end{itemize}
\end{thm}

The proof is based on a technique of higher order expansion  \cite{W. E(1995),C. Wang(2000)}. It is  very complex  and   postponed to the Appendix.

The following result is on  the convergence of damped Newton's iteration \eqref{equ:numnum1}-\eqref{Newton}.
\begin{thm}
Suppose the initial data $f_0(X)\in\mathcal{E}_{M}$ is   positive and bounded for $X\in\mathcal{Q}$,  then Newton's iteration \eqref{equ:numnum1}-\eqref{Newton} is convergent in $\mathcal{Q}$.
\end{thm}
We can first prove that  $J(y)$, defined in \eqref{DisEnergy}, is a self-concordant function  \cite{C.H. Duan(2017), Y. Nesterov(1994)}.
Then based on  Theorem\ 2.2.3 in \cite{Y. Nesterov(1994)}, damped Newton's iteration \eqref{equ:numnum1}-\eqref{Newton} is convergent in $ \mathcal{Q}$.  We omit the details.

%

%

\section{Numerical Results}
\label{sec:numerical results}
In this section, we show some numerical results. To demonstrate the accuracy of the numerical schemes, in the first example,  we solve  a problem with a smooth solution. In the second example, we consider a free boundary problem with a exact Barenblatt solution. We check the convergence  for the solution and the finite speed of propagation. In the third example, we focus on numerical simulation for the waiting time. Finally we report some results for problems with two support sets at the initial state in Example 4.

The error of a numerical solution  is measured  in the $\mathcal{L}^2$ and $\mathcal{L}^{\infty}$ norms defined as:
\begin{equation}\label{L2}
\|e_h\|_2^2=\frac{1}{2}\left(e_{h_0}^2 h_{x_{0}}+\sum\limits_{i=1}^{M-1}e_{h_i}^2 h_{x_i}+e_{h_M}^2 h_{x_M}\right),
\end{equation}
and
\begin{equation}\label{Linf}
\|e_h\|_{\infty}=\max\limits_{0\leq i\leq M}\{|e_{h_i}|\},
\end{equation}
 where $e_h=(e_{h_0},e_{h_1},\cdots,e_{h_M})$ and
for the error of the density  $f-f_h$,
 \begin{equation*}
 h_{x_i}=x_{i+1}-x_{i-1}, \ \ 1\leq i \leq M-1; \ \ \
 h_{x_0}=x_{1}-x_{0}; \ \
 h_{x_M}=x_{M}-x_{M-1}, \end{equation*}
and for  the error of the  trajectory $x-x_h$,
 \begin{equation}
h_{x_i}=2h,  \ \ 1\leq i\leq M-1, \ \ \ h_{x_0}=h_{x_M}=h.\nonumber
\end{equation}
 where $h$ is the spatial step.

\noindent{\bf Example 1.} \label{Positive} Convergence rate for problem with smooth solution

Consider the problem \eqref{eqPMEori}-\eqref{eqPMEboun} in dimension one with a smooth positive initial state
\begin{eqnarray}
\label{ini1}
&& f_0(x)=\sin(\pi x)+0.5, \ \ x\in\Omega=(0,1).
\end{eqnarray}
We solve the trajectory equation \eqref{eqtra2} in Case 1 (\eqref{eqtra3} in Case 2) with the initial and boundary condition \eqref{eqtraboun}-\eqref{eqtraini} by the fully discrete scheme \eqref{eqtranum2} in Case 1 (\eqref{eqtranum3} in Case 2) and approximate the density function $f$  in \eqref{equ:conservationL} by \eqref{numdist}. The reference 'exact' solution is obtained numerically on a much fine mesh with $h=\frac{1}{100000},\    \tau=\frac{1}{100000}$.

Tables \ref{tablePME11} and  \ref{tablePMEL11}  show  the convergence rate in Cases 1 and 2, respectively.  The rate for  density $f$ and trajectory $x$ in the   $\mathcal{L}^2$ and  $\mathcal{L}^{\infty}$ norm is  2nd order in space and 1st order in time for each scheme. But  the linear scheme \eqref{eqtranum3} in Case 2 is more efficient. 

\begin{threeparttable}[b]
\scriptsize
\centering
\setlength{\abovecaptionskip}{10pt}
\setlength{\belowcaptionskip}{-3pt}
\caption{\scriptsize {\bf Example 1}.  Convergence  rate of solution $f$ and trajectory $x$ in Case 1 at final time $T=0.05$}
\label{tablePME11}
\begin{tabular}{@{ } l c c c c c c c c c c}

\hline
\multicolumn{1}{l}{}
&\multicolumn{8}{c}{$m=\frac{5}{3}$}\\\hline
 $M$    &$\tau$ &$ \mathcal{L}^2$-error $ (f) $    & Order   &$ \mathcal{L}^{\infty}$-error$ (f) $    & Order &$ \mathcal{L}^2$-error $(x)$   & Order  &$ \mathcal{L}^{\infty}$-error $(x)$   & Order & CPU (s)  \\\hline
100 &1/100 & 1.1304e-02&  &1.6847e-02& &1.5122e-03 & &2.2356e-03& &0.1872\\\hline
 200 &1/400 &2.6730e-03 & 2.1144&3.8606e-03&2.1820&3.5665e-04 &2.1200&5.2869e-04&2.1143&0.6084\\\hline
 400 &1/1600 &6.4528e-04 &2.0712 &9.2707e-04&2.0821&8.6042e-05 &2.0725&1.2761e-04&2.0716&2.1840\\\hline
 800 &1/6400 &1.5246e-04 &2.1163 &2.1878e-04&2.1187&2.0324e-05 &2.1167&3.0145e-05&2.1165&8.7361\\\hline

\hline
\multicolumn{1}{l}{}
&\multicolumn{8}{c}{$m=2$}\\\hline
  $M$    &$\tau$ &$ \mathcal{L}^2$-error $ (f) $    & Order   &$ \mathcal{L}^{\infty}$-error$ (f) $    & Order &$ \mathcal{L}^2$-error $(x)$   & Order  &$ \mathcal{L}^{\infty}$-error $(x)$   & Order & CPU (s)  \\\hline
100 &1/100 & 8.4443e-03&  &1.2463e-02& &1.1269e-03 & &1.1269e-03& &0.1716\\\hline
 200 &1/400 &1.8021e-03 & 2.3429&2.5826e-03&2.4129&2.3982e-04 &2.3494&2.3982e-04&2.3494&0.5304\\\hline
 400 &1/1600 &4.1921e-04 &2.1495 &5.9831e-04&2.1583&5.5749e-05 &2.1509&5.5749e-05&2.1509&2.0748\\\hline
 800 &1/6400 &9.8039e-05 &2.1379 &1.3980e-04&2.1399&1.3034e-05 &2.1386&1.3034e-05&2.1386&8.0185\\\hline
\end{tabular}
\begin{tablenotes}
     \scriptsize
        \item[1]   $ \mathcal{L}^2$-error and $ \mathcal{L}^{\infty}$-error  is defined by  \eqref{L2} and \eqref{Linf}, respectively.

                \item[2]    $\tau$ is the time step and  $h=\frac{1}{M}$ is the space step.
        \item[3]  CPU (s) is the CPU time (seconds).
        \end{tablenotes}
\end{threeparttable}

 \begin{threeparttable}[b]
\scriptsize
\centering
\setlength{\abovecaptionskip}{10pt}
\setlength{\belowcaptionskip}{-3pt}
\caption{\scriptsize {\bf Example 1}. Convergence  rate of  solution $f$ and trajectory $x$  in Case 2 at final time $T=0.05$}

\label{tablePMEL11}
\begin{tabular}{@{ } l c c c c c c c c c c}

\hline
\multicolumn{1}{l}{}
&\multicolumn{8}{c}{$m=\frac{5}{3}$}\\\hline
 $M$    &$\tau$ &$ \mathcal{L}^2$-error  $ (f) $    &Order   &$ \mathcal{L}^{\infty}$-error $ (f) $    & Order &$ \mathcal{L}^2$-error $(x)$   &Order    &$ \mathcal{L}^{\infty}$-error $(x)$   &Order  & CPU(s) \\\hline
100 &1/100 & 1.0617e-02&       &1.6396e-02&       &1.4212e-03 &      &2.0955e-03&  &0.0000    \\\hline
 200 &1/400 &2.5002e-03 & 2.1233&3.6535e-03&2.2439&3.3374e-04 &2.1291&4.9444e-04&2.1190&0.0000 \\\hline
 400 &1/1600 &6.0295e-04 &2.0733 &8.7321e-04&2.0920&8.0425e-05 &2.0749&1.1922e-04&2.0736&1.5600e-02\\\hline
 800 &1/6400 &1.4238e-04 &2.1174 &2.0580e-04&2.1215&1.8987e-05 &2.1179&2.8150e-05&2.1176&6.2400e-02\\\hline

\hline
\multicolumn{1}{l}{}
&\multicolumn{8}{c}{$m=2$}\\\hline
  $M$    &$\tau$ &$ \mathcal{L}^2$-error  $ (f) $    &Order   &$ \mathcal{L}^{\infty}$-error $ (f) $    & Order &$ \mathcal{L}^2$-error $(x)$   &Order    &$ \mathcal{L}^{\infty}$-error $(x)$   &Order  & CPU(s) \\\hline
100 &1/100 & 8.0516e-03&        &1.2168e-02&      &1.0750e-03 &      &1.5887e-03& &0.0000\\\hline
 200 &1/400 &1.7134e-03 & 2.3497&2.4675e-03&2.4656&2.2803e-04 &2.3572&3.3833e-04&2.3479&0.0000\\\hline
 400 &1/1600 &3.9861e-04 &2.1492 &5.7051e-04&2.1625&5.3010e-05 &2.1508&7.8690e-05&2.1498&1.5600e-02\\\hline
 800 &1/6400 &9.3216e-05 &2.1381 &1.3324e-04&2.1409&1.2392e-05 &2.1388&1.8397e-05&2.1386&6.2400e-02\\\hline
\end{tabular}

\begin{tablenotes}
     \scriptsize
        \item[1]   $ \mathcal{L}^2$-error and $ \mathcal{L}^{\infty}$-error  is defined by  \eqref{L2} and \eqref{Linf}, respectively.

                \item[2]    $\tau$ is the time step and  $h=\frac{1}{M}$ is the space step.
        \item[3]  CPU (s) is the CPU time (seconds).
        \end{tablenotes}
\end{threeparttable}

\noindent{\bf Example 2} \label{ex:BS}   Numerical  finite propagation speed  for problem with free boundary

Barenblatt solution \cite{G.I. Barenblatt (1952), R.E. Pattle (1959),J. L. Vazquez(2007),Ya.B. Zeldovich (1950)} in dimension one can be expressed by
  \begin{equation}
     B_m(x,t)=(t+1)^{-k}\Big(1-\frac{k(m-1)}{2m}\frac{|x|^2}{(t+1)^{2k}}\Big)^{1/(m-1)}_{+}, \  x\in \mathbb{R},\ t\leq0,
     \label{equ:BSE}
  \end{equation}
  where $l_{+}=\max\{l,0\}$ and $k=(m+1)^{-1}$. The solution has a compact
support $[-\xi^B_m(t), \xi^B_m(t)]\subsetneqq\Omega$ with the interface $|x|=\xi^B_m(t)$ moving outward in a finite speed,
where
\begin{equation}
\xi^B_m(t):=\sqrt{\frac{2m}{k(m-1)}}\cdot (t+1)^k.
\label{equ:BSInterE}
\end{equation}

Let the computing domain be $\Omega  = (-10, 10)$. We take Barenblatt profile  $B_m(x,0)$   as the initial data in problem \eqref{eqPMEori}-\eqref{eqPMEboun}. For a finite time interval,  the interface can not reach the boundary of $\Omega$, so the boundary condition \eqref{eqPMEboun} is valid.   We solve the trajectory equation \eqref{eqtraC2}-\eqref{eqiniC2} in Case 1 (\eqref{eqtraC3}-\eqref{eqiniC3} in Case 2) by the fully discrete scheme \eqref{equtrafreenum2}-\eqref{equtrafreenumboun2} in Case 1(\eqref{equtrafreenum3}-\eqref{equtrafreenumboun3} in Case 2).

\begin{figure}
\captionsetup{font={scriptsize}}
\centering
\subfigure[\scriptsize $t=2$]
{\includegraphics[width=6cm,height=5cm]{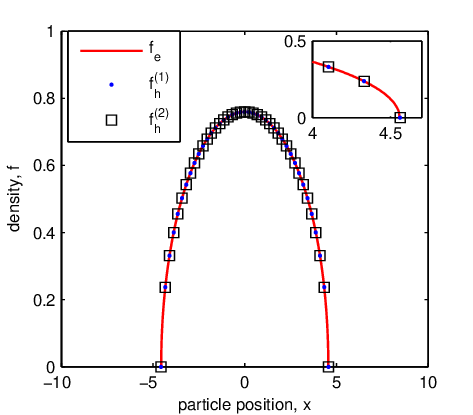}}\hspace{.1in}
\subfigure[\scriptsize $t=10$]
{\includegraphics[width=6cm,height=5cm]{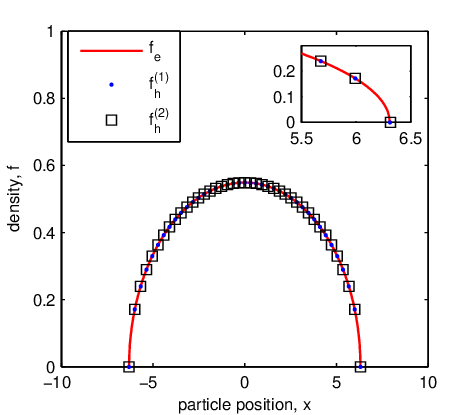}}\hspace{.1in}
\caption{\textbf{Example 2}. The evolution of $f$; $f_e$ is the exact solution; $f_h^{(1)}$ and $f_h^{(2)}$ are numerical  solutions in Case 1 and Case 2, respectively ($m=3$, $M=2000$, $\tau=1/1000$)}
\label{fig:PMELBSDen}
\end{figure}
\begin{figure}
\captionsetup{font={scriptsize}}
\centering
\subfigure[\scriptsize  The trajectory in Case 1]
{\includegraphics[width=6cm,height=5cm]{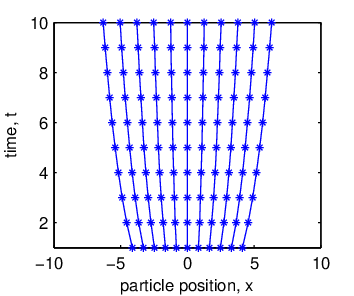}}\hspace{.1in}
\subfigure[\scriptsize  The trajectory in Case 2]
{\includegraphics[width=6cm,height=5cm]{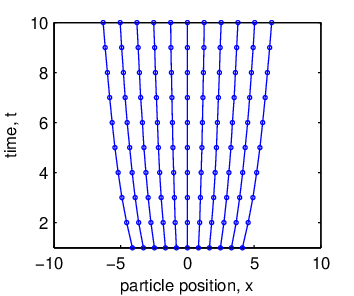}}\hspace{.1in}
\caption{\textbf{Example 2}.  The evolution of particle position for  $m=3$ over time  ($M=2000$, $\tau=1/1000$)}
\label{fig:PMEBSParticle}
\end{figure}
\begin{figure}
\captionsetup{font={scriptsize}}
\centering
\includegraphics[width=8cm,height=7cm]{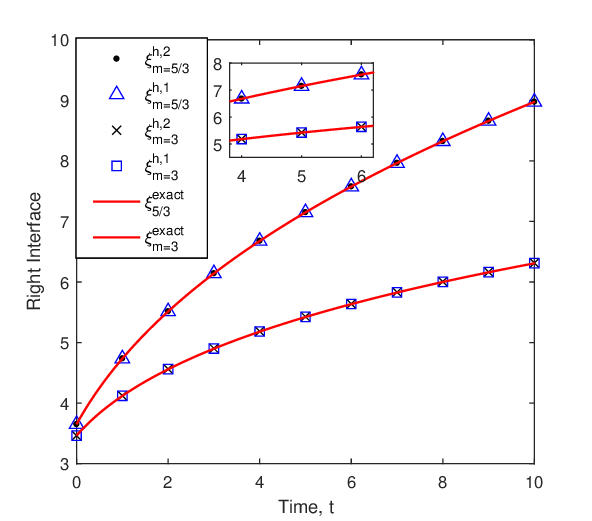}
\caption{\textbf{Example 2}.  The evolution of the right interface over time for  different $m$ ($M=2000$ and   $\tau=1/1000$);  $\xi^{h,1}$ and $\xi^{h,2}$ denote the numerical interfaces in Case 1 and Case 2, respectively ($M=2000$, $\tau=1/1000$)}
\label{fig:PMELBSInterface}
\end{figure}

  Fig. \ref{fig:PMELBSDen} shows the  numerical    and  exact solutions for $m=3$ at time $t=2$ and $t=10$. The results demonstrate that  the numerical solutions in Case 1 and Case 2 can approximate to the exact solution without oscillation.   The evolution of  the trajectory in  both cases   over time  for $m=3$ is  shown in Fig. \ref{fig:PMEBSParticle}:  particles move outward in a finite speed without twisting or exchanging.  Fig. \ref{fig:PMELBSInterface} shows the evolution of  the right interface  for numerical solutions and the exact solution with different $m$ ($m=\frac{5}{3}$,  $m=3$) in Case 1 and Case 2. Table \ref{tableBSInterface} shows the error of the  right interface with different $m$ ($m=\frac{5}{3}$, $m=2$,  $m=3$, $m=5$) at time $T=1$. The results mean that the numerical interface in each case is a  good approximation to the exact one and moves in a finite speed.
   %

  Table \ref{tableBS}  shows the convergence rate of $f$ in Case 1 and Case 2. We present
  the  error   in  $\mathcal{L}^2$  norm and the error at $X=0$  at  time $T=1$  for  $m=5/3$ and $m=3$.
The results show that  the convergence rate is deteriorated when $m$ is getting large. This is due to the deteriorated regularity of the solution.    The error of $f$ at $X=0$ keeps the rate of  2nd order since $f$ is still smooth far away from the interface. Both numerical schemes have the same rate,  but the  error of $f$ in Case 2 is larger. Table \ref{tableBSNum} shows the convergence rate of $f$ in $\mathcal{L}^{\infty}$ norm for $m=\frac{5}{3}$ in the three cases: the  numerical schemes  lead to the same convergence rate i.e., $1$st order.

\begin{threeparttable}[b]
\scriptsize
\setlength{\abovecaptionskip}{10pt}
\setlength{\belowcaptionskip}{-3pt}
\caption{\scriptsize \textbf{Example 2}. The error of right interface $\xi_r$  at  $T=1$}
\label{tableBSInterface}
\begin{tabular}{@{ } l c c c c}

\hline
 $m$ &$\frac{5}{3}$ &2&3 &5 \\\hline
 $|\xi^{h,1}_r-\xi^{exact}_r|$ & 6.6911e-04 & 2.3153e-04 & 2.9872e-03&5.1647e-03\\\hline
 $|\xi^{h,2}_r-\xi^{exact}_r|$& 9.6066e-04&3.9205e-03&6.2808e-03 &6.8532e-03 \\\hline
\end{tabular}
\begin{tablenotes}
       \scriptsize
        \item[1] $\xi^{exact}_r$  denotes the exact right interface; $\xi^{h,1}_r$ and $\xi^{h,2}_r$ denote the numerical right interfaces in Case 1 and Case 2, respectively.
        \end{tablenotes}
\end{threeparttable}

\begin{threeparttable}[b]
\tiny
\centering
\setlength{\abovecaptionskip}{10pt}
\setlength{\belowcaptionskip}{-5pt}
\setlength{\abovecaptionskip}{0.cm}
\setlength{\belowcaptionskip}{0.cm}

\caption{\scriptsize {\bf Example 2.} The convergence rate of  $f$  at the finite time $T=1$}
\label{tableBS}
\begin{tabular}{@{ } l c c c c c c c c c}

\hline
\multicolumn{1}{l}{}
&\multicolumn{7}{c}{$m=\frac{5}{3}$}\\\hline
 $M$    &$\tau$  &$\mathcal{L}^2$-error $\big(f_h^{(1)}\big)$    &Order  &$\mathcal{L}^2$-error $\big(f_h^{(2)}\big)$     & Order    & Error at $X=0$ $\big(f_h^{(1)}\big)$ &Order & Error at $X=0$ $\big(f_h^{(2)}\big)$ &Order 
 \\\hline
1000 &1/250 &5.6454e-05 & &6.1225e-04 &      &2.5417e-05 & &2.7701e-04 & 
\\\hline
2000 &1/1000 &1.4133e-05 &1.9972 &1.5281e-04 &2.0033 &6.3626e-06&1.9974&6.9154e-05&2.0029 
\\\hline
4000 &1/4000 &3.5351e-06 &1.9990 &3.8184e-05 &2.0009  &1.5912e-06&1.9993 &1.7282e-05&2.0007  
\\\hline
 8000 &1/16000 &8.8404e-07 &1.9994 &9.5445e-06 &2.0003 &3.9782e-07&1.9998&4.3202e-06&2.0002  
 \\\hline

 \multicolumn{1}{l}{}
&\multicolumn{7}{c}{$m=3$}\\\hline
$M$    &$\tau$  &$\mathcal{L}^2$-error $\big(f_h^{(1)}\big)$      &Order  &$\mathcal{L}^2$-error $\big(f_h^{(2)}\big)$     & Order    & Error at $X=0$ $\big(f_h^{(1)}\big)$ &Order & Error at $X=0$ $\big(f_h^{(2)}\big)$ &Order\\\hline
1000 &1/250 &1.3480e-03 &               &5.8979e-03 &            &4.1682e-05 & &1.9361e-04 &  
\\\hline
 2000 &1/1000 &6.7614e-04 &0.9969 &2.4952e-03 &1.1819 &1.0821e-05&1.9259&4.9069e-05&1.9728 
 \\\hline
 4000 &1/4000 &3.4194e-04 &0.9887 &1.2050e-03 &1.0353 &2.8617e-06&1.8907 &1.2570e-05&1.9519 
 \\\hline
 8000 &1/16000 &1.7310e-04 &0.9877 &6.0168e-04 &1.0014 &7.7215e-07&1.8531&3.2253e-06&1.9306 
 \\\hline

\end{tabular}
\begin{tablenotes}
       \scriptsize
       \item[1]    $f_h^{(1)}$ and $f_h^{(2)}$ are the numerical solutions of the problem \eqref{eqPMEori}-\eqref{eqPMEboun} in Case 1 and Case 2, respectively.
        \item[2]   $ \mathcal{L}^2$-error $ \big(f_h^{(i)}\big)$ is the error of $f_h^{(i)}$ in $\mathcal{L}^2$ norm  defined by \eqref{L2},   $i=1,2$.  
        \item[3]    $\tau$ is the time step; $h=\frac{1}{M}$ is  the space step.

        \end{tablenotes}

\end{threeparttable}

\begin{threeparttable}[b]
 \scriptsize
\centering
\setlength{\abovecaptionskip}{10pt}
\setlength{\belowcaptionskip}{-3pt}
\caption{\scriptsize {\bf Example 2.} The convergent rate of $f$ in $\mathcal{L}^{\infty}$  norm at final time $T=1$}

\label{tableBSNum}
\begin{tabular}{@{ } l c c c c c c c }
\hline
\multicolumn{1}{l}{}
&\multicolumn{7}{c}{$m=\frac{5}{3}$}\\\hline
 $M$    &$\tau$  &$\mathcal{L}^{\infty}$-error $ \big(f_h^{(0)}\big)$       &Order  &$\mathcal{L}^{\infty}$-error  $\big(f_h^{(1)}\big)$     & Order   &$\mathcal{L}^{\infty}$-error  $\big(f_h^{(2)}\big)$     & Order  
 \\\hline
100 &1/10 &3.44e-04 &       &1.00e-03 &       &7.46e-03 & 
\\\hline
250 &1/25 &9.82e-05 &1.37 &2.65e-04 &1.45 &2.83e-03&1.06
\\\hline
1000 &1/100 &1.32e-05 &1.45 &6.28e-05 &1.04 &6.93e-04&1.01 
\\\hline
2500 &1/250 &3.40e-06 &1.48 &2.51e-05 &1.00 &2.76e-04&1.00 
 \\\hline

\end{tabular}
\begin{tablenotes}
       \scriptsize
       \item[1] $f_h^{(0)}$  is the numerical solution by VPS \cite{M. Westdickenberg(2010)};  $f_h^{(1)}$ and $f_h^{(2)}$ are the numerical solutions in Case 1 and Case 2, respectively.
        \item[2]   $ \mathcal{L}^{\infty}$-error $ (f_i) $, $i=0,1,2$ are the error of  solution $f$ in $\mathcal{L}^{\infty}$  norm defined by \eqref{Linf}.
        \end{tablenotes}

\end{threeparttable}

\noindent{\bf Example 3}\label{ex:WT} Numerical simulation for the waiting time

The waiting-time phenomenon  occurs for a certain type of initial states \cite{J. L. Vazquez(2007)}.  Without loss of generality we consider  the left interface. Similar  argument can be obtained for  the right interface.   Recalling the  trajectory equation  \eqref{eqbouC2} or  \eqref{eqbouC3} at the left interface, we have
\begin{equation}\label{equ:freeBou}
 \partial_t x=-\frac{m}{m-1}\frac{\partial_X[f_0(X)^{m-1}]}{(\partial_X x)^m}, \mbox{at}\  X=\xi_1^0,  \ t>0,
\end{equation}
where   $f_0(X)$ is the smooth initial state with compact support $[\xi_1^0, \xi_2^0]$. At the initial time, $\partial_X x\equiv 1$, so if $\partial_X[f_0(X)^{m-1}] = 0$ at $X=\xi_1^0$, then $x_t(\xi_1^0,0) = 0$ and it is possible to have a positive waiting time.

 If the left interface keeps waiting till time $t^*>0$, then $\xi_1^t \equiv \xi_1^0$, for $t\le t^*$. This means that we must have, at $X=\xi_1^0$, $\partial_t x \equiv  0$, for  $t < t^*$ and $\partial_t x < 0$, for $t=t^*+\epsilon$ with any sufficiently small $\epsilon >0$.
Hence  the waiting time can be characterized as:
\begin{equation}\label{WTC}
t^*:=\inf\Big\{t>0: x_t= - \frac{m}{m-1}\frac{\partial_X[f_0(X)^{m-1}]}{(\partial_X x)^m} < 0,\  \mbox{as}\  X\rightarrow \xi_1^0. \Big\}.
\end{equation}
Noting that, at $X= \xi_1^0$, the numerator $\partial_X[f_0(X)^{m-1}]$ is fixed and only the denominator   $(\partial_X x)^m$ changes when time evolves. If there exists a positive waiting time $t^* >0$, we must have that, at $X= \xi_1^0$,  $\partial_X[f_0(X)^{m-1}]=0$  and as time evolves, $(\partial_X x)^m$ becomes smaller and smaller and comes to the same order infinitesimal as $\partial_X[f_0(X)^{m-1}]$ as $X\rightarrow \xi_1^0$ at time $t=t^*$. So we have another criterion for the waiting time:
\begin{equation} \label{WTn}  t^* \mbox{\ is\ the\ first\ time\ instant\ when\ }  \mathcal{B}(t):= \frac{\partial_X[f_0(X)^{m-1}]}{(\partial_X x)^m}\ \mbox{is\ not\ a\ infinitesimal\ as\ } X\rightarrow \xi_1^0.
 \end{equation}
%
%
%
Next we focus on finding the  criterion for the numerical waiting time $t^*_h$.  Let $$\mathcal{B}_h^n:=\frac{\bar{D}_{h}[(f_0(X_0))^{m-1}]}{(\bar{D}_{h} x^{n}_{h,0})^m},$$
where  the difference operator $\bar{D}_{h}$ is defined in \eqref{deri1} and  $x^{n}_h=(x_{h,0}^{n},\cdots,x_{h,M}^{n})$ is the numerical trajectory  position at time $t^{n}$, $n=0,\cdots,N$.

The numerical waiting time $t^*_h$   is determined by  the following  criterion:
 \begin{equation}\label{WTC2_1}
t^*_h:=\min\Big\{t^n: \Big|\frac{\mathcal{B}^n_{2h}}{\mathcal{B}_h^n}\Big|\leq 1 \Big\}.
\end{equation}
To get $\mathcal{B}_{2h}^n$ in the above formula, we need to know the trajectory $x^n_{2h}$. we don't need to solve the trajectory problem again by spacial step $2h$. We just select it from the given solution $x^n_h$, i.e., $x^n_{2h}=(x^n_{h,0}, x^n_{h,2}, x^n_{h,4},\cdots, )$.
\begin{rem}
The numerical criterion \eqref{WTC2_1} is an approximation of the continuous criterion \eqref{WTn} in the sense that if $\mathcal{B}(t)$ is  infinitesimal  as $X\rightarrow \xi_1^0$, then $\mathcal{B}_h< \mathcal{B}_{2h}$ for any  sufficiently small $h>0$.
\end{rem}

Now we present the algorithm for problems with waiting time.

\noindent{\bf Algorithm for Waiting time}
\begin{itemize}
\item {\bf Step 1.}\ \ For time $t^n, n=0,1, \cdots$, solve the trajectory equation \eqref{eqtraC2}-\eqref{eqiniC2}  in Case 1 (\eqref{eqtraC3}-\eqref{eqiniC3} in Case 2)  by the fully discrete scheme  \eqref{equtrafreenum2}-\eqref{equtrafreenumboun2} in Case 1 (\eqref{equtrafreenum3}-\eqref{equtrafreenumboun3} in Case 2) but replacing the boundary condition  \eqref{eqbouC2}  in Case 1 (\eqref{eqbouC3} in Case 2) by $\partial_t x = 0$ and replacing the boundary condition \eqref{equtrafreenumboun2}  in Case 1 (\eqref{equtrafreenumboun3} in Case 2) by $x^{n+1}_0 = \xi_1^0, \ x^{n+1}_M = \xi_2^0$.

    Check the criterion \eqref{WTC2_1} for $x^{n+1}$. If it is not valid, goto next time step. If it is valid, then set $t^*_h = t^{n+1}$. $n^*= n+1$ and  goto Step 2.
\item {\bf Step 2.}\ \ For time $t^n, n=n^*,n^*+1, \cdots$, solve the trajectory equation  \eqref{eqtraC2}-\eqref{eqiniC2} in Case 1 (\eqref{eqtraC3}-\eqref{eqiniC3} in Case 2)   by the fully discrete scheme \eqref{equtrafreenum2}-\eqref{equtrafreenumboun2} in Case 1 (\eqref{equtrafreenum3}-\eqref{equtrafreenumboun3} in Case 2) .
\end{itemize}

Now we consider the following data set-up:
\begin{eqnarray}
&&
 \Omega=(-5,5),\\
&&
f_0(x)=\left\{
\begin{array}{lcl}
\big\{\frac{m-1}{m}[(1-\theta)\sin^2(x)+\theta\sin^4(x)\big\}^{1/(m-1)}, &&{x\in[-\pi, 0]},\\
0, &&{\mbox{otherwise\ in} \  \Omega },
\end{array}\right.
\label{WT}
\end{eqnarray}
where $\theta\in[0,\frac{1}{4}]$. Then the waiting time is positive  and
the exact one  \cite{D. G. Aronson(1983)} is:
\begin{equation}\label{equ:ExWT}
t^*_{exact}=\frac{1}{2(m+1)(1-\theta)}.
\end{equation}

\begin{figure}
\captionsetup{font={scriptsize}}
\centering
\subfigure[\scriptsize $t=0$]
{\includegraphics[width=4cm,height=3cm]{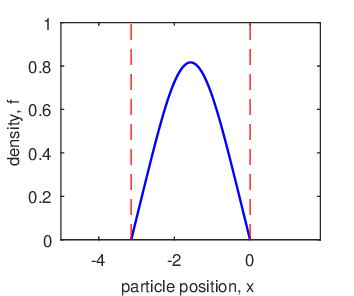}}\hspace{.3in}
\subfigure[\scriptsize $t=0.05$]
{\includegraphics[width=4cm,height=3cm]{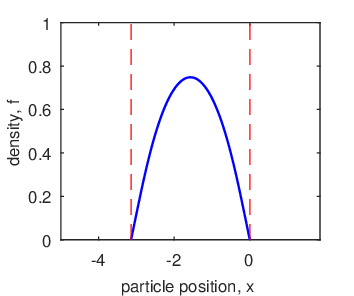}}\hspace{.3in}
\subfigure[\scriptsize $t=0.1$]
{\includegraphics[width=4cm,height=3cm]{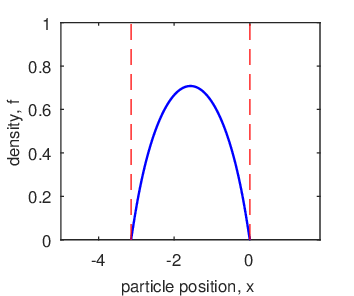}}\hspace{1in}\\
\subfigure[\scriptsize $t=0.169$]
{\includegraphics[width=4cm,height=3cm]{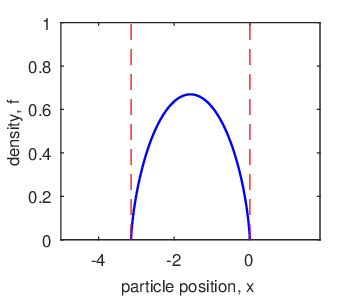}}\hspace{.3in}
\subfigure[\scriptsize $t=0.5$]
{\includegraphics[width=4cm,height=3cm]{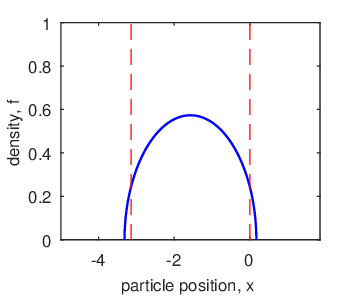}}\hspace{.3in}
\subfigure[\scriptsize $t=1$]
{\includegraphics[width=4cm,height=3cm]{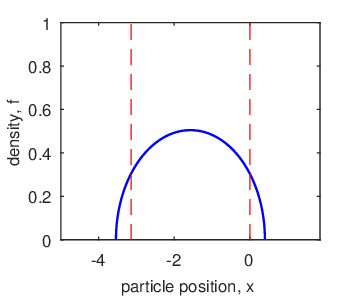}}\hspace{1in}
\caption{ {\bf Example 3. Waiting time:} Evolution of solution $f$ in Case 1 ($m=3$, $M=1000$,  $\tau=1/2000$)}
\label{fig:WTDen}
\end{figure}

\begin{figure}
\captionsetup{font={scriptsize}}
\centering
\subfigure[\scriptsize   the influence of  $\theta$ ($m=3$)]
{\includegraphics[width=6cm,height=5cm]{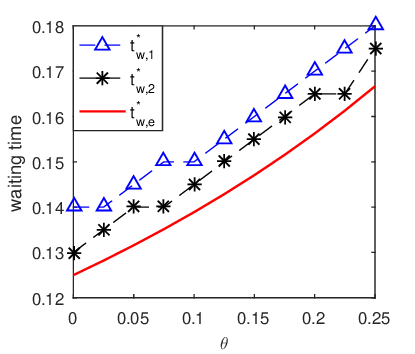}}\hspace{.1in}
\subfigure[\scriptsize  the influence of   $m$  ($\theta=\frac{1}{4}$)]
{\includegraphics[width=6cm,height=5cm]{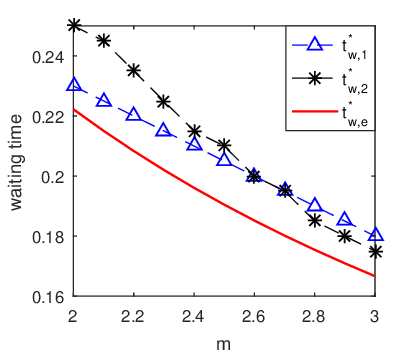}}\hspace{.1in}
\caption{ {\bf Example 3. Waiting time:}  the influence of  $\theta$ and  $m$;  $t^*_e$ is the exact waiting time given by  \eqref{equ:ExWT}; $t^*_{w,1}$ and $t^*_{w,2}$  are the numerical waiting time in Case 1 and Case 2, respectively ($M=200$, $\tau=1/200$)}
\label{fig:PMELWT}
\end{figure}

Fig. \ref{fig:WTDen} depicts that the evolution of numerical solution $f$   over grid with spatial step $h=\pi/M$ ($M=1000$) and the time step $\tau=1/2000$ for $m=3$ and $\theta=\frac{1}{4}$ in Case 1. The results show that the  waiting time does exist. After the time about $0.169$, the interface moves outward in a finite speed. In the whole process, we obtain the numerical solution without oscillation.   Fig. \ref{fig:PMELWT} (a) and (b) present the comparison of the numerical and exact waiting time for different $\theta$ and  $m$ in Case 1 and Case 2.  The results show that the  numerical waiting time    is a  good approximation to the exact one in each case. 
 Furthermore, Table \ref{tableWTL} presents the  error of waiting time for $m=3$ and $\theta=\frac{1}{4}$ over different grids ($M=500$, $\tau=1/1000$; $M=1000$, $\tau=1/2000$; $M=2000$, $\tau=1/4000$; $M=4000$, $\tau=1/8000$) in Case 1 and Case 2. It  shows that the numerical waiting time is convergent to the exact one in each case.

\begin{threeparttable}[b]
\scriptsize
\centering
\setlength{\abovecaptionskip}{10pt}
\setlength{\belowcaptionskip}{-3pt}
\caption{\scriptsize \textbf{Example 3}.  The convergence rate of  waiting time ($m=3$, $\theta=\frac{1}{4}$)}
\label{tableWTL}
\begin{tabular}{@{ } l c c c c c c c c c}

\hline
 $M$ &$\tau$ &$t^*_{w,1}$ &$|t^*_{w,1}-t^*_{w,e}|$ &Order& $CPU^{1}$(s)&$t^*_{w,2}$&$|t^*_{w,2}-t^*_{w,e}|$ &Order& $CPU^{2}$(s)\\\hline
 25 & $\frac{1}{25}$  &0.24&0.0733&  &4.6875e-02& 0.24 & 0.0733&  &1.5625e-02\\\hline
 50& 1/50 &0.20&0.0333&1.1006&4.6875e-02  &0.20&0.0333 &1.1006&3.1250e-02  \\\hline
100&1/100&0.19&0.0233&0.7146 &7.8125e-02 & 0.18&0.0133&1.2519&3.1250e-02\\\hline
200&1/200&0.180&0.0133&0.8759 &1.4063e-01& 0.175&0.0083&0.8012&4.6875e-02 \\\hline
$t^*_{w,e}$ & & 0.16667 & & &0.16667 \\\hline
\end{tabular}
\begin{tablenotes}
       \scriptsize
        \item[1]  $t^*_{w,e}$ is the exact waiting time by \eqref{equ:ExWT}; $t^*_{w,1}$  and $t^*_{w,2}$ are the waiting time in Case 1 and 2, respectively.
                  \item[2] $CPU^{1}$(s) and  $CPU^{2}$(s)  denote the CPU time (seconds)  in Case 1 and 2, respectively.
                \end{tablenotes}
\end{threeparttable}

\noindent{\bf Example 4}\label{ex:twoSupport} Numerical simulation for problem with two separate support sets at initial time

  We now consider a problem with a step function as the initial  state. In   problem \eqref{eqPMEori}-\eqref{eqPMEboun},  let $m=5$,  $\Omega=(-5,5)$ and
\begin{equation}
f_0(x)=\left\{
\begin{array}{lcl}
1, &&{x\in(0.5, 3)},\\
1.5, &&{x\in(-3,-0.5)},\\
0, && \mbox{otherwise}.
\end{array}\right.
\label{CC}
\end{equation}

The example models the movement and interaction of two supports. Before the two supports meet, we solve two problems independently. When the two supports meet at time $t^*_m$, we should reconstruct the two parts of solution into a whole with single support over an equidistance mesh and then take it as initial state to solve problem \eqref{eqtraC2}-\eqref{eqiniC2}  in Case 1 (\eqref{eqtraC3}-\eqref{eqiniC3} in Case 2)  starting from $t=t^*_m$.

\begin{figure}
\captionsetup{font={scriptsize}}
\centering
\subfigure[\scriptsize $t=0$]
{\includegraphics[width=4cm,height=3cm]{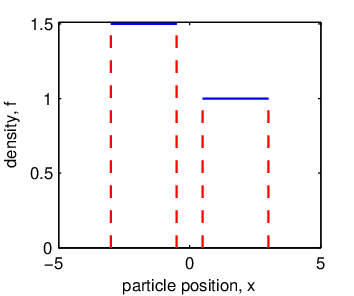}}\hspace{1in}
\subfigure[\scriptsize $t=0.1$]
{\includegraphics[width=4cm,height=3cm]{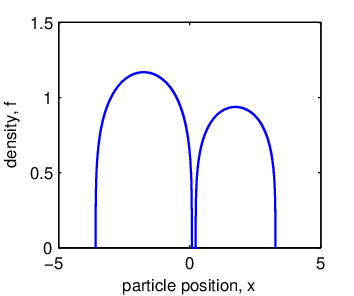}}\hspace{1in}\\
\subfigure[\scriptsize meet at $t=0.1415$]
{\includegraphics[width=4cm,height=3cm]{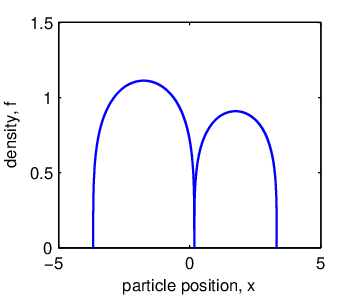}}\hspace{1in}
\subfigure[\scriptsize reconstruction at $t=0.1415$]
{\includegraphics[width=4cm,height=3cm]{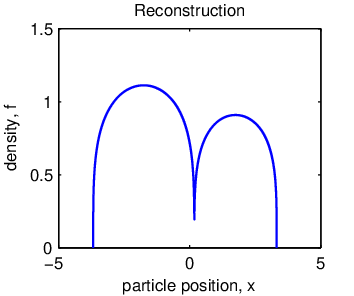}}\hspace{1in}
\subfigure[\scriptsize $t=0.5$]
{\includegraphics[width=4cm,height=3cm]{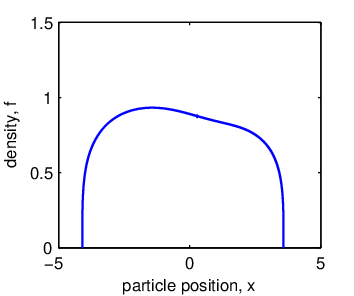}}\hspace{1in}
\subfigure[\scriptsize $t=1$]
{\includegraphics[width=4cm,height=3cm]{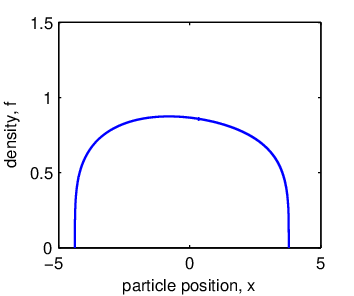}}\hspace{1in}\\
\subfigure[\scriptsize $t=2$]
{\includegraphics[width=4cm,height=3cm]{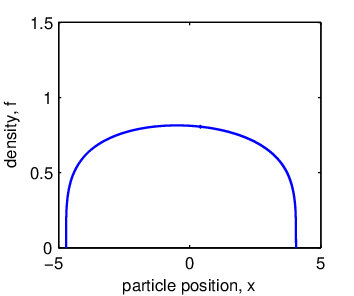}}\hspace{1in}
\subfigure[\scriptsize particle position]
{\includegraphics[width=4cm,height=3cm]{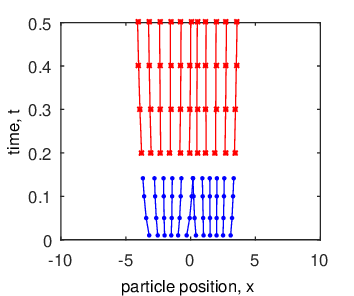}}\hspace{1in}\\
\caption{ Example 4. The evolution of density $f$ and particle position $x$ in Case 1 for  $m=5$  over time  ($M=5000$,  $\tau=1/10000$ and the space grid size of  reconstruction $M_2=10000$)}
\label{fig:TwoDenDH}
\end{figure}

The spatial step is chosen as $h=(3-0.5)/M$ ($M=5000$) for each support and the time step is $\tau=1/10000$. In   Case 1,   Figs.\ref{fig:TwoDenDH} (a)$\sim$(c)  show that as time evolves, the two supports  expand and   meet at time $t^*_m=0.1415$.  At this time,  a reconstruction is taken  by monotone piecewise cubic interpolation \cite{F. N. Fritsch(1980)} over an equidistance grid  with partition number $M_2=10000$, shown in  Fig.\ref{fig:TwoDenDH}(d).  Figs.\ref{fig:TwoDenDH}(e)$\sim$(g) show the evolution after meeting.
Oscillations do not appear around the free boundary during the whole process.  Fig. \ref{fig:TwoDenDH}(h) shows the movement  of particles  in this process. The numerical solution in  Case 2 has the  similar results and the meeting time is $t^*_m=0.1383$.
\begin{rem} \label{rem:monotone}
 The meeting time of two supports is defined as:
\begin{equation}
t^*_m:=\inf\limits_{t>0}\{|x^l_M-x^r_0|\leq 10^{-10}\},
\end{equation}
where $x^l_M$ is the endpoint of the left support and $x^r_0$ the first point of right support.

\end{rem}

\section{Concluding Remarks}
   In this paper, the numerical solution of the PME based on EnVarA has been proposed and analyzed. Originated from the different energy dissipation laws of the system, we mainly consider two numerical schemes of  the trajectory equation obtained by the balance of LAP and MDP. Based on the total energy density $f\ln f$, the proposed numerical scheme is proven to be uniquely solvable on an admissible convex set, mainly thanks to the singularity of the total energy. Based on the  total energy density $\frac{1}{2f}$, the numerical scheme is linear. In turn, the energy dissipation rate of both schemes has been an outcome of the variational approach. Moreover, the second order convergence in space and the first order convergence in time have been theoretically justified  for both schemes, with a careful application of higher order asymptotic expansion of the numerical schemes to obtain higher order consistency. 
   According to the numerical simulation results for both schemes, no oscillation appears around the free boundary, and  the  finite propagation speed could be numerically calculated. We also give a criterion that can compute the waiting time and numerical convergence of the waiting time is reported,  which is the first such result for this problem. Furthermore, the numerical scheme based on $\frac{1}{2f}$ is linear and  more efficient.

One obvious limitation of this work is associated with the one-dimensional nature of the problem. Solving for multi-dimensional PME by this energetic method will be left to our future works.

\noindent\emph{Acknowledgments.}
This work is supported in part by NSF of China under the grants 11271281. Chun Liu and Cheng Wang are partially supported by NSF grants DMS-1216938, DMS-1418689, respectively.

\section*{Appendix: Proof of Theorem \ref{convergence}}
 \renewcommand{\theequation}{A.\arabic{equation}}

  \renewcommand{\thethm}{A.\arabic{thm}}

We first prove the Theorem \ref{convergence} for the numerical scheme  \eqref{eqtranum2}.
Before that,  we introduce a higher order  approximate expansion  of the exact solution since a regular  expansion (a second order in space and a first order in time) do not obtain the convergence rate in Theorem \ref{convergence}.

\begin{lem}
Assume a higher order  approximate solution  of the exact solution $x_e$:
\begin{equation}\label{approximation}
W:=x_e+\tau w_\tau^{(1)}+\tau^2 w^{(2)}_\tau +h^2 w_h,
\end{equation}
where $w_\tau^{(1)}$, $w_\tau^{(2)}$, $w_h \in C^{\infty}{(\Omega; 0,T)}$. Then there exists a small $\tau_0>0$, such that $\forall\tau, h \leq \tau_0$, $\widetilde{D}_h W>0$, i.e., $W\in\mathcal{Q}$, where $\tau$ and  $h$ are the time step and the spatial step, respectively.
\end{lem}
\noindent{\textbf{Proof:}}
Since a point-wise level of $x_e\in\mathcal{Q}$, i.e.,  $\exists$ $\varepsilon_0>0$, such that $D_h x_e>\varepsilon_0>0$.
For small $\tau_0$, such that $\|\tau D_h w_\tau^{(1)}\|_{L^{\infty}}\leq\frac{1}{9}\varepsilon_0$, $\|\tau^2 D_h w_\tau^{(2)}\|_{L^{\infty}}\leq\frac{1}{9}\varepsilon_0$ and $\|h^2 D_h w_h\|_{L^{\infty}}\leq\frac{1}{9}\varepsilon_0$, for $\forall\tau, h \leq \tau_0$.
As a consequence, for $\forall \tau, h\leq\tau_0$, we have
\begin{equation}\label{U_Q}
D_h W\geq\frac{1}{3}\varepsilon_0>0,
\end{equation}
then $W\in\mathcal{Q}$.   $\hfill\Box$

Then we proceed into the proof of Theorem \ref{convergence}.

\noindent\textbf{Proof of Theorem \ref{convergence}}:
A careful Taylor expansion with high order of  \eqref{eqtra2} in both time and space  shows that
\begin{eqnarray}
&&
 \frac{f_{0}(X_{i})}{m(\frac{f_0(X_{i})}{\widetilde{D}_h x_{e_{i}}^{n}})^{m-1}}\frac{x_{e_{i}}^{n+1}-x_{e_{i}}^{n}}{\tau}=-d_{h}\Big(\frac{f_{0}(X)}{D_{h} x_{e}^{n+1}}\Big)_{i}+ \tau l^{(1)}_i +\tau^{2}l^{(2)}_i+\tau^{3}l^{(3)}_i+h^{2}g^{(1)}_i+h^4g^{(2)}_i,\nonumber\\
&&
\ \ \ \ \ 1\leq i\leq M-1,\nonumber\\
&&
 \quad \mbox{with $x_{e_{0}}^{n+1} = X_0$ ,   $\quad  x_{e_{M}}^{n+1} = X_M $}, \label{consistency-2}
\end{eqnarray}
where $\|l^{(1)}\|_2$, $\|l^{(2)}\|_2$, $\|l^{(3)}\|_2$, $\|g^{(1)}\|_2$, $\|g^{(2)}\|_2\leq C_e$,  with $C_e$  dependent on the exact solution.

To perform a higher order consistency analysis for an approximate solution  of the exact solution, we have to construct the approximation $W$ as in \eqref{approximation}.

The term $w_\tau^{(1)}\in C^{\infty}{(\Omega; 0,T)}$ is given by the following linear equation:
\begin{eqnarray}
  && \frac{f_0(X)}{m(\frac{f_0(X)}{\partial_X x_e})^{m-1}}\partial_t w_\tau^{(1)}+  \frac{m-1}{m (\frac{f_0(X)}{\partial_X x_e})^{m-2}} \partial_t x_e\cdot\partial_X w_\tau^{(1)}
  =\partial_X\Big(\frac{f_0(X)}{(\partial_X x_e)^2}\partial_X w_\tau^{(1)}\Big)-l^{(1)}, \label{equ:u_1}\notag\\
  &&
 w_\tau^{(1)}|_{\partial\Omega}=0,  \label{u_1_bou} \ \ \ \   w_\tau^{(1)}(\cdot,0)=0.
  \label{u_1_ini}
\end{eqnarray}

The term $w^{(2)}_\tau\in C^{\infty}{(\Omega; 0,T)}$ is given by the following linear equation:
\begin{eqnarray}
  && \frac{f_0(X)}{m(\frac{f_0(X)}{\partial_X x_e})^{m-1}}\partial_t w^{(2)}_\tau+  \frac{m-1}{m (\frac{f_0(X)}{\partial_X x_e})^{m-2}} \partial_t x_e\cdot\partial_X w^{(2)}_\tau \nonumber\\
  &&
  +\frac{(m-1)(m-2)}{2m (\frac{f_0(X)}{\partial_X x_e})^{m-2}\partial_X x_e}(\partial_X w^{(1)}_\tau)^2\cdot\partial_t x_e +\frac{(m-1)}{m (\frac{f_0(X)}{\partial_X x_e})^{m-2}}\partial_t w^{(1)}_\tau\cdot \partial_X w^{(1)}_\tau\nonumber \\
  &&
=\partial_X\Big(\frac{f_0(X)}{(\partial_X x_e)^2}\partial_X w^{(2)}_\tau\Big)-\partial_X\Big(\frac{f_0(X)}{(\partial_X x_e)^3}(\partial_X w^{(1)}_\tau)^2\Big)-l^{(2)},\label{equ:u_2}\notag \\
  &&
  w^{(2)}_\tau|_{\partial\Omega}=0,  \label{u_2_bou} \ \ \ \  w^{(2)}_\tau(\cdot,0)=0.
  \label{u_2_ini}
\end{eqnarray}

The term $w_h\in C^{\infty}{(\Omega; 0,T)}$ is given by the following linear equation:
\begin{eqnarray}
  &&
  \frac{f_0(X)}{m(\frac{f_0(X)}{\partial_X x_e})^{m-1}}\partial_t w_h+  \frac{(m-1)\partial_t x_e}{m (\frac{f_0(X)}{\partial_X x_e})^{m-2}} \partial_X w_h=\partial_X\Big(\frac{f_0(X)}{(\partial_X x_e)^2}\partial_X w_h\Big)-g^{(1)}, \label{equ:u_h}\notag \\
  &&
  w_h|_{\partial\Omega}=0, \label{u_h_bou}\ \ \ \   w_h(\cdot,0)=0.
  \label{u_h_ini}
\end{eqnarray}
Since  $w_\tau^{(1)}$, $w_\tau^{(2)}$, $w_h$ are only depend on $W$ and $x_e$, we have the following estimate:
\begin{eqnarray}
  &&
  \|W-x_e\|_{H^m}=\tau\|w_\tau^{(1)}\|_{H^m}+\tau^2\|w_\tau^{(1)}\|_{H^m}+h^2\|w_h\|_{H^m}\leq C'(\tau+h^2). \label{estimate_u1}
  %
\end{eqnarray}

With such an expansion term, the constructed approximation $W\in\mathcal{Q}$ satisfies the numerical scheme with a higher order truncation error:
\begin{eqnarray}\label{num:U}
&&
\frac{f_{0}(X_i)}{m\big(\frac{f_0(X)}{\widetilde{D}_h W^{n}}\big)_i^{m-1}}\cdot\frac{W^{n+1}_i-W^{n}_i}{\tau}=-d_{h}\Big(\frac{f_{0}(X)}{D_{h} W^{n+1}}\Big)_i+ \tau^{3}l^{*}_i+h^4g^{*}_i, \ \ 1\leq i\leq M-1, \nonumber\\
&&
\mbox{with\ \ } W_0^{n+1}=X_0, \ \ W_M^{n+1}=X_M, \ \ n=0,1,\cdots, N-1,
\end{eqnarray}
where $l^{*}$, $g^{*}$ are dependent on $l^{(1)}$, $l^{(2)}$, $l^{(3)}$, $g^{(1)}$, $g^{(2)}$ and the derivatives of $w^{(1)}_\tau$, $w^{(2)}_\tau$, $w_h$.

Then we define $\tilde{e}_i^n := W_i^n - x_{h_i}^n$, $0 \le i \le M$, $n=0,1,\cdots, N$. In other words, instead of a direct comparison between the numerical solution and exact PDE solution, we evaluate the numerical error between the numerical solution and the constructed solution $W$. The higher order truncation error enables us to obtain a required $W_h^{1,\infty}$ of the numerical solution, which is necessary in the nonlinear convergence analysis.

Note that the discrete $L^2$ norm $\|\tilde{e}^{0}\|_{2}=0$ at time step $t^0$. We assume at time step $t^n$:
\begin{equation}\label{priori}
\|\tilde{e}^{n}\|_2\leq \gamma(\tau^3+h^4),
\end{equation}
where   the constant $\gamma$, given in \eqref{M}, is dependent on the exact solution $x_e$ and its derivative.

Then we have the following estimates:
\begin{eqnarray}
&&
\|\widetilde{D}_h\tilde{e}^{n}\|_{2}\leq \gamma(\tau^2+h^3), \label{a-priori}\\
&&
 \|\widetilde{D}_h \tilde{e}^{n} \|_{\infty} \leq C_m\frac{\|\widetilde{D}_h\tilde{e}^{n}\|_{2}}{h^{1/2}}\leq C_m \gamma(\tau^{\frac{3}{2}}+h^{\frac{5}{2}}), \ \mbox{if\ \ } h=O(\tau),
\label{a-priori_h} \\
&&
\|\widetilde{D}_h x^{n}_h \|_{\infty}=\|\widetilde{D}_h W^n-\widetilde{D}_h\tilde{e}^n\|_{\infty}\leq C^*+1 :=C_0^*, \label{a-priori_x}\\
&&
\mbox{with\ \ } C^*:=\|\widetilde{D}_h W^n\|_{\infty},\ \ \mbox{if\ \ } C_m\gamma(\tau^{\frac{3}{2}}+h^{\frac{5}{2}})\leq 1, \nonumber\\
&&
\Big\|\frac{\widetilde{D}_h x^{n}_h-\widetilde{D}_h x^{n-1}_h }{\tau}\Big\|_{\infty}=\Big\|\frac{\widetilde{D}_h W^{n}-\widetilde{D}_h W^{n-1} }{\tau}-\frac{\widetilde{D}_h \tilde{e}^{n}-\widetilde{D}_h \tilde{e}^{n-1} }{\tau}\Big\|_{\infty}\leq \tilde{C}^*_t+1,\label{a-priori_x_t}\\
&&
\mbox{with\ \ } \tilde{C}^*_t:=\Big\|\frac{\widetilde{D}_h W^{n}-\widetilde{D}_h W^{n-1} }{\tau}\Big\|_{\infty}, \ \mbox{if\ \ } C_m\gamma(\tau^{\frac{1}{2}}+h^{\frac{3}{2}})\leq 1. \label{a-priori-diff-time}
\end{eqnarray}
%
For $x_h,W\in\mathcal{Q}$, i.e., $\exists$ $\delta_0>0$, such that  $\widetilde{D}_h W^n_i\geq \delta_0$, then  $\widetilde{D}_h x^n_{h_i}\geq \frac{\delta_0}{2}>0$, $0\leq i\leq M$, if  $C_m \gamma(\tau^{\frac{3}{2}}+h^{\frac{5}{2}})\leq\frac{\delta_0}{2}$.



In turn, subtracting \eqref{num:U} from the numerical scheme  \eqref{eqtranum2} yields
\begin{eqnarray}
  && \ \ \ \frac{f_0(X_i)}{m\big(\frac{f_0(X)}{\widetilde{D}_h x^{n}_{h}}\big)_i^{m-1}}\cdot\frac{\tilde{e}^{n+1}_i - \tilde{e}^n_i}{\tau} +
  \frac{f_0(X_i)}{ m[f_0(X_i)]^{m-1}}\cdot\frac{W^{n+1}_i-W^{n}_i}{\tau}\cdot[\big(\widetilde{D}_h W^{n}\big)_i^{m-1}-\big(\widetilde{D}_h x_h^{n}\big)_i^{m-1}]
  \nonumber \\
  &&
  = d_h \left( \frac{f_0(X)}{D_h W_i^{n+1} D_h x_h^{n+1}}  D_h \tilde{e}^{n+1}\right)_i   + \tau^3 f^*_i+h^4g^{*}_i,\ \ 1\leq i\leq M-1,  \nonumber \\
  && \mbox{with\ }\tilde{e}^{n+1}_0 = \tilde{e}^{n+1}_M = 0,
  \label{consistency-3}
\end{eqnarray}
in which the form of the  left term comes from the following identity:
 \begin{eqnarray}
 &&
 \ \ \ \frac{f_{0}(X_i)}{m\big(\frac{f_0(X)}{\widetilde{D}_h W^{n}}\big)_i^{m-1}}\frac{W^{n+1}_i-W^{n}_i}{\tau}-
 \frac{f_{0}(X)}{m\big(\frac{f_0(X_i)}{\widetilde{D}_h x_h^{n}}\big)_i^{m-1}}\frac{x^{n+1}_{h_i}-x^{n}_{h_i}}{\tau}
  \nonumber\\
 &&
 = \frac{f_{0}(X_i)}{\tau m [f_0(X_i)]^{m-1}}[(\widetilde{D}_h W^{n})_i^{m-1}(W^{n+1}_i-W^{n}_i)
 -(\widetilde{D}_h x^{n}_{h})_i^{m-1}(x^{n+1}_{h_i}-x^{n}_{h_i}) \nonumber\\
 &&
 \  \  \ +\big(\widetilde{D}_h x^{n}_h\big)_i^{m-1}(W^{n+1}_i-W^{n}_i)-\big(\widetilde{D}_h x^{n}_h\big)_i^{m-1}(W^{n+1}_i-W^{n}_i)]\nonumber\\
 &&
 =\frac{f_0(X_i)}{ m[f_0(X_i)]^{m-1}}\cdot\frac{W^{n+1}_i-W^{n}_i}{\tau}\cdot[\big(\widetilde{D}_h W^{n}\big)_i^{m-1}-\big(\widetilde{D}_h x_h^{n}\big)_i^{m-1}] \nonumber \\
 &&
 \ \ \ +\frac{f_0(X_i)}{m\big(\frac{f_0(X)}{\widetilde{D}_h x_h^{n}}\big)_i^{m-1}}  \cdot\frac{\tilde{e}_i^{n+1} - \tilde{e}_i^n}{\tau}. \nonumber
  \end{eqnarray}

Based on the preliminary results about the discrete inner-product, taking a discrete inner product with (\ref{consistency-3}) by $2 \tilde{e}^{n+1}$ gives
\begin{eqnarray}
&&
 2\left\langle \alpha_n (\tilde{e}^{n+1} - \tilde{e}^n), \tilde{e}^{n+1}\right\rangle -2\tau\left\langle  d_h \left( \frac{f_0(X)}{D_h W^{n+1} D_h x^{n+1}_h}  D_h \tilde{e}^{n+1}\right), \tilde{e}^{n+1}\right\rangle\nonumber\\
 &&
 =-2\tau\left\langle \frac{f_0(X)}{ m[f_0(X)]^{m-1}}\cdot\frac{W^{n+1}-W^{n}}{\tau}\cdot[(\widetilde{D}_h W^{n})^{m-1}-(\widetilde{D}_h x^{n}_h)^{m-1}], \tilde{e}^{n+1}\right\rangle \nonumber\\
 &&
 \ \ \ +2\tau\left\langle \tau^3 f^*+h^4g^{*},\tilde{e}^{n+1}\right\rangle,
  \label{convergence-1}
\end{eqnarray}
where  \begin{equation}\label{alpha}\alpha_n:=\frac{f_0(x)}{m\big(\frac{f_0(X)}{\widetilde{D}_h x^{n}_h}\big)^{m-1}}.\end{equation}

For the first  term of left side, we have
\begin{equation}
\begin{split}
2\left\langle \alpha_n  (\tilde{e}^{n+1} - \tilde{e}^n),\tilde{e}^{n+1}\right\rangle
&
=\alpha_n \|\tilde{e}^{n+1}\|^2_{2}+\alpha_n \|\tilde{e}^{n+1}-\tilde{e}^{n}\|^2_{2}-\alpha_n \|\tilde{e}^{n}\|^2_{2} \\
&
\geq \alpha_n \|\tilde{e}^{n+1}\|^2_{2}-\alpha_n \|\tilde{e}^{n}\|^2_{2}. \label{convergence-proof-1}
\end{split}
\end{equation}

For the second term of left side, we have
\begin{eqnarray}
&&
 -2\tau\left\langle  d_h \left( \frac{f_0(X)}{D_h W^{n+1} D_h x^{n+1}_h}  D_h \tilde{e}^{n+1}\right), \tilde{e}^{n+1}\right\rangle \nonumber\\
&&
=2\tau\left\langle \frac{f_0(X)}{D_h W^{n+1} D_h x_h^{n+1} }  D_h \tilde{e}^{n+1}, D_h \tilde{e}^{n+1}\right\rangle_e \geq 0,
%
\label{convergence-proof-2} 
\end{eqnarray}
in which  the summation by parts formula (\ref{FD-inner product-3-1}) is applied with $\tilde{e}^{n+1}_0 = \tilde{e}^{n+1}_N = 0$.

For the right side term, we have
\begin{eqnarray}
&&
-2\tau\left\langle \frac{f_0(X)}{ m[f_0(X)]^{m-1}}\cdot\frac{W^{n+1}-W^{n}}{\tau}\cdot[(\widetilde{D}_h W^{n})^{m-1}-(\widetilde{D}_h x^{n}_h)^{m-1}], \tilde{e}^{n+1}\right\rangle \nonumber\\
&&
=-2\tau\left\langle \frac{f_0(X)}{m[f_0(X)]^{m-1}}\cdot\frac{W^{n+1}-W^{n}}{\tau}\cdot[(m-1)(\widetilde{D}_h \zeta^{n})^{m-2}\widetilde{D}_h \tilde{e}^{n}], \tilde{e}^{n+1}\right\rangle \nonumber\\
&&
\leq 2\tau C_1\|\widetilde{D}_h  \tilde{e}^{n}\|_2\|\tilde{e}^{n+1}\|_2, \ \ \Big(C_1:=\frac{(m-1)B_fC^*_tC_{\zeta}^{m-2}}{mb_f^{(m-1)}} \Big),\\
&&
\leq\tau C_1\|\widetilde{D}_h  \tilde{e}^{n}\|_2^2+\tau C_1\|\tilde{e}^{n+1}\|_2^2, \nonumber
 \label{convergence-proof-3}
\end{eqnarray}
in which $ C^*_t=\|W_t\|_{\infty}$,
$\widetilde{D}_h\zeta$ is between $\widetilde{D}_h x_h^n$ and $\widetilde{D}_h W^n$ and $\|\widetilde{D}_h\zeta\|_{\infty}\leq C_\zeta$ with
\begin{equation*}
C_\zeta:=
\begin{cases}
C^*_0, & m\geq2,\nonumber\\
\delta_0/2, & m<2.
\end{cases}
\end{equation*}

The local truncation error term could be bounded by the standard Caught inequality:
\begin{equation}
\begin{split}
2\tau\left\langle \tau^3 f^*+h^4g^{*}, \tilde{e}^{n+1}\right\rangle
&
\leq \tau\|\tau^3 f^*+h^4g^{*}\|_2^2+\tau\|\tilde{e}^{n+1}\|^2_2  \\
&
\leq \tau C(\tau^3+h^4)^2+\tau\|\tilde{e}^{n+1}\|^2_2.
  \label{convergence-proof-4}
  \end{split}
\end{equation}



Next we estimate $\|D_h x^{n+1}_h\|_{\infty}$ roughly.
Based on \eqref{alpha}, $\alpha_n$ can be estimated by
$$C_\alpha:=\frac{b_f}{m(\frac{B_f}{\delta_0/2})^{m-1}}\leq\|\alpha\|_n\leq\frac{B_f}{mb_f^{m-1}}(C_0^*)^{m-1}:=\bar{C}_\alpha.$$
A substitution of \eqref{convergence-proof-1} - \eqref{convergence-proof-4} into (\ref{convergence-1}) with \eqref{a-priori} leads to
\begin{equation*}
\begin{split}
 (\alpha_n-\tau (1+C_1))\| \tilde{e}^{n+1} \|_{2}^2
 &
 \leq \alpha_n\| \tilde{e}^n \|_{2}^2+\tau C_1 \|\widetilde{D}_h  \tilde{e}^{n}\|_2^2+\tau C(\tau^3+h^4)^2 \nonumber\\
 &
 \leq \tau \bar{C}(\tau^2+h^3)^2,
  \end{split}
\end{equation*}
where $\bar{C}$ is dependent on $C_1$, $\gamma$ and $\bar{C}_\alpha$.
Then we have
\begin{eqnarray}
&&
  \| \tilde{e}^{n+1} \|_{2}^2 \leq \widetilde{C}^2\tau(\tau^2+h^3)^2, 
  \mbox{\ \ i.e.,\ \ } \| \tilde{e}^{n+1} \|_{2} \leq \widetilde{C}\tau^{\frac{1}{2}}(\tau^2+h^3), \end{eqnarray}
  \mbox{with\ \ } $\widetilde{C}:=\Big(\frac{\bar{C}}{C_\alpha/2}\Big)^{\frac{1}{2}}$, \ \ \
 \mbox{if\ \ } $\tau (1+C_1)\leq \bar{C}_\alpha/2$.

Based on  the inverse inequality \eqref{equ:inverse}, we obtain that, by choosing  $h=O(\tau)$,
\begin{equation}
\| \tilde{e}^{n+1} \|_{\infty}\leq \frac{C_m\| \tilde{e}^{n+1} \|_{2}}{h^{\frac{1}{2}}}\leq C_m\widetilde{C}(\tau^2+h^3).
\end{equation}
Then we have
\begin{eqnarray}
\|D_h {x}^{n+1}_h \|_{\infty}=\|D_h {W}^{n+1} -D_h \tilde{e}^{n+1} \|_{\infty}\leq C^*+C_m\widetilde{C}(\tau^2+h^3)\leq C^*+1:=C^*_0, \label{rough_x}\ \ \
\end{eqnarray}
if $C_m\widetilde{C}(\tau^2+h^3)\leq 1$.

Then \eqref{convergence-proof-2} can be re-estimated as following:
\begin{equation}\label{convergence-proof-2-2}
2\tau\left\langle \frac{f_0(X)}{D_h W^{n+1} D_h x_h^{n+1}}  D_h \tilde{e}^{n+1}, D_h \tilde{e}^{n+1}\right\rangle_e
\geq 2\tau C_2\|D_h \tilde{e}^{n+1}\|_2^2,
\end{equation}
with $C_2:=\frac{b_f}{C^*C_0^*}$.

As a consequence, a substitution of \eqref{convergence-proof-1} - \eqref{convergence-proof-4} with \eqref{convergence-proof-2-2} into (\ref{convergence-1}) leads to
\begin{eqnarray}
\alpha_n\| \tilde{e}^{n+1} \|_{2}^2 - \alpha_n\| \tilde{e}^n \|_{2}^2   +
 \tau C_2\|\widetilde{D}_h \tilde{e}^{n}\|^2_{2} 
 \leq \tau\Big(1+\frac{C_1^2}{C_2}\Big)\|\tilde{e}^{n+1}\|^2_{2}
 +\tau C (\tau^3+h^4)^2,\nonumber
\end{eqnarray}
where the following estimates are applied: $\|\widetilde{D}_h x^n_h\|_2\leq\|D_h x^n_h\|_2$ and
\begin{equation}
2\tau C_1\|\widetilde{D}_h  \tilde{e}^{n}\|_2\|\tilde{e}^{n+1}\|_2\leq
\tau\frac{C_1^2}{C_2}||\tilde{e}^{n+1}||^2_{2} +\tau C_2||\widetilde{D}_h \tilde{e}^{n}||^2_{2},
\end{equation}

Then summing in time shows that 
\begin{eqnarray}
 \alpha_n \| \tilde{e}^{n+1} \|_{2}^2+
  \tau C_2\sum\limits_{k=1}^{n+1}||\widetilde{D}_h \tilde{e}^{k}||^2_{2}
 & \leq &\tau\sum\limits_{k=1}^{n}\frac{(\alpha_k- \alpha_{k-1})}{\tau}\|\tilde{e}^{k}\|^2_2+
  \tau \Big(\frac{C_1^2}{C_2}+1\Big)\sum\limits_{k=1}^{n+1}\|\tilde{e}^{k}\|^2_{2} \notag\\
  && + CT (\tau^3+h^4)^2. \notag \\
\| \tilde{e}^{n+1} \|_{2}^2+
  \tau \frac{C_2}{C_\alpha}\sum\limits_{k=1}^{n+1}||\widetilde{D}_h \tilde{e}^{k}||^2_{2}
& \leq&  \frac{\tau}{C_\alpha} (\frac{C_1^2}{C_2}+1+\widetilde{C}_{\alpha})\sum\limits_{k=1}^{n+1}\|\tilde{e}^{k}\|^2_{2} +\frac{CT}{C_\alpha} (\tau^3+h^4)^2,\notag
\end{eqnarray}
where we have used the estimate
\begin{equation}
\begin{split}
\big\|\frac{\alpha^k-\alpha^{k-1}}{\tau}\big\|_{\infty}
&
=\Big\|\frac{f_0(X)}{m[f_0(X)]^{m-1}}\cdot\frac{(\widetilde{D}_h x^{k}_h)^{m-1}-(\widetilde{D}_h x^{k-1}_h)^{m-1}}{\tau}\Big\|_{\infty} \nonumber\\
&
=\Big\|\frac{f_0(X)}{m[f_0(X)]^{m-1}}(m-1)(\widetilde{D}_h \vartheta)^{m-2}\frac{\widetilde{D}_h x^{k}_h-\widetilde{D}_h x^{k-1}_h}{\tau}\Big\|_{\infty}\nonumber\\
&
\leq \frac{(m-1)B_f}{m b_f^{m-1}}(C_{\vartheta})^{m-2}(\tilde{C}^*_t+1):=\widetilde{C}_{\alpha},\nonumber
 \end{split}
\end{equation}
in which $T$ is the terminal time, \eqref{a-priori-diff-time} is applied and  $\widetilde{D}_h \vartheta$ is between $\widetilde{D}_h x^k_h$ and $\widetilde{D}_h x^{k-1}_h$ with
\begin{equation*}
\|\widetilde{D}_h \vartheta\|_{\infty}\leq C_{\vartheta}:=
\begin{cases}
C^*_0, & m\geq2,\\
\delta_0/2, & m<2.
\end{cases}
\end{equation*}

In turn, an application of discrete Gronwall inequality yields the desired convergence result:
\begin{eqnarray}
&& \| \tilde{e}^{n+1} \|_{2}^2   +
  \tau \frac{C_2}{C_\alpha}\sum\limits_{k=1}^{n+1}||\widetilde{D}_h \tilde{e}^{k}||^2_{2}
  \leq e^{TC_0}\frac{CT}{C_\alpha}(\tau^3+h^4)^2, \nonumber\\
  && i.e.,  \| \tilde{e}^{n+1} \|_{2}\leq \gamma(\tau^3+h^4),\nonumber
\end{eqnarray}
where  $C_0:=\frac{1}{C_\alpha} (\frac{C_1^2}{C_2}+\widetilde{C}_{\alpha}+1)$ and
\begin{equation}
\gamma:=\Big(\frac{CT }{C_\alpha}\Big)^{\frac{1}{2}}e^{\frac{C_0T}{2}}.\label{M}
\end{equation}
Therefore a-priori assumption \eqref{a-priori} is valid at $t^{n+1}$:
\begin{equation}
\| \tilde{e}^{n+1} \|_{2}\leq \gamma(\tau^3+h^4),
\end{equation}
if $\tau\leq\min\big\{(\frac{1}{\gamma C_m})^{2},\frac{\bar{C}_{\alpha}}{2(1+C_1)}, (\frac{1}{C_m\widetilde{C}})^{\frac{1}{2}},(\frac{\delta_0}{2C_m\gamma})^{\frac{2}{3}}\big\}$.

 Based on
  \begin{equation}\label{derivative}
 \|\widetilde{D}_h \tilde{e}^{n+1}\|_2 =\|\widetilde{D}_h x^{n+1}_h-\widetilde{D}_h W^{n+1}\|_2\leq \gamma(\tau^2+h^3),
 \end{equation}
 we obtain
  \begin{equation}\label{derivative_u_e}
  \|\widetilde{D}_h x_h^{n+1}-  \widetilde{D}_h x_e^{n+1}\|_2\leq C(\tau+h^2).
  \end{equation}
Next we focus on  the error  between the numerical solution $f_h^{n+1}$ and the exact solution $f_e^{n+1}$ of the problem \eqref{eqPMEori}-\eqref{eqPMEboun}.
 \begin{equation}
 \begin{split}
\|f_{e}^{n+1}-f_{h}^{n+1}\|_2&=\Big\|\frac{f_0(X)}{\partial_X x^{n+1}_e}-\frac{f_0(X)}{\widetilde{D}_h x^{n+1}_h}\Big\|_2 \nonumber \\
&
=\left\|\frac{f_0(X)}{\partial_X x^{n+1}_e}-\frac{f_0(X)}{\widetilde{D}_h x^{n+1}_e}+\frac{f_0(X)}{\widetilde{D}_h x^{n+1}_e}-\frac{f_0(X)}{\widetilde{D}_h x^{n+1}_h}\right\|_2  \nonumber \\
&
\leq C(\tau+h^2).
\end{split}
\end{equation}

With some minor modifications, we can prove Theorem \ref{convergence} for the numerical scheme  \eqref{eqtranum3}. Due to the linear scheme  \eqref{eqtranum3},  the rough estimation  of  $\|D_h x^{n+1}_h\|_{\infty}$ can  not  be used.
 $\hfill\Box$

{\footnotesize
\bibliographystyle{unsrt}

\begin{thebibliography}{0}
\bibitem{D.G. Aronson(1969)}
D.G. Aronson, Regularity properties of flows through porous media, SIAM J. Appl. Math. 17 (1969) 461-467.




\bibitem{D. G. Aronson(1983)}
D.G. Aronson, L.A. Caffarelli,  S. Kamin, How an initially stationary interface begins to move in porous medium flow, SIAM J. Math. Anal. 14 (4)  (1983) 639-658.



\bibitem{G.I. Barenblatt (1952)}
G.I. Barenblatt, On some unsteady motions of a liquid or a gas in a porous medium, Prikl. Mat. Mekh. 16 (1) (1952) 67-78 (in Russian).

\bibitem{M. Bertsch (1990)}
M. Bertsch,  R. Dal Passo, A numerical treatment of a super degenerate equation with applications to the porous media equation, Quart. Appl. Math.  48 (1990) 133-152.

\bibitem{E. DiBenedetto(1984)}
E. DiBenedetto, D. Hoff,  An interface tracking algorithm for the porous medium equation, Trans. Am. Math. Soc. 284 (1984)  463-500.




\bibitem{Q. Du(2009)}
Q. Du,  C. Liu, R. Ryham,  X. Wang, Energetic variational approaches in modeling vesicle and fluid interactions, Physica. D. 238 (2009) 923-930.

\bibitem{C.H. Duan(2017)}
C. Duan, C. Liu, C. Wang,  X.  Yue, Numerical complete solution for random genetic drift by Energetic Variational approach, arXiv:1803.09436 (2018).

\bibitem{W. E(1995)}
W. E, J. -G. Liu, Projection method I: convergence and numerical boundary layers, SIAM J. Numer. Anal. 32 (1995) 1017-1057.




\bibitem{B. Eisenberg(2010)}
B. Eisenberg, Y.K. Hyon,  C. Liu, Energy variational analysis of ions in water and channels: Field theory for primitive models of complex ionic fluids,  J. Chem. Phys. 133 (10)  (2010) 104.

\bibitem{F. N. Fritsch(1980)}
 F.N. Fritsch, R.E. Carlson, Monotone Piecewise Cubic Interpolation, SIAM J. Numer. Anal. 17 (1980) 238-246.

 \bibitem{J. L. Graveleau(1971)}
J.L. Graveleau,  P. Jamet, A finite difference approach to some degenerate nonlinear parabolic equations, SIAM J.  Appl.  Math. 20 (1971) 199-223.

\bibitem{J. Gratton(1998)}
J. Gratton, C. Vigo, Evolution of self-similarity, and other properties of waiting-time solutions of the porous medium equation: the case of viscous gravity currents, J. Appl. Math. 9 (1998) 327-350.

\bibitem{T. Huang(2016)}
T. Huang, F. Lin, C. Liu, and C. Wang,  Finite time singularity of the nematic liquid crystal flow in dimension three, Arch. Ration. Mech. An.  221(3) (2016) 1223-1254.
%

\bibitem{Y. Hyon(2010)}
Y. Hyon, D.Y. Kwak and C. Liu, Energetic variational approach in complex fluids: maximum dissipation principle, Discrete Contin. Dyn. Syst.  26(4) (2010) 1291-1304

\bibitem{S. Jin(1998)}
S. Jin, L. Pareschi,  G. Toscani, Diffusive relaxation schemes for multi-scale discrete-velocity kinetic equations,  SIAM  J. Numer. Anal. 35  (6) (1998) 2405-2439.

\bibitem{A.S. Kalasnikov(1967)}
A.S. Kala\v{s}nikov, Formation of singularities in solutions of the equation of nonstationary filtration, \v{Z}. Vy\v{c}isl. Mat. Mat. Fiz. 7 (1967) 440-444.

\bibitem{Hajime Koba(2017)}
 H. Koba, C. Liu, Y. Giga, Energetic variational approaches for incompressible fluid systems on an evolving surface, Quart. Appl. Math. 75 (2017)  359-389.


\bibitem{L. S. Leibenzon(1930)}
L.S. Leibenzon, The motion of a gas in a porous medium, Complete works, vol 2, Acad.
Sciences URSS, Moscow, Russian (1953).

\bibitem{C. Liu(2003)}
C. Liu, J. Shen, A phase field model for the mixture of two incompressible fluids and its approximation by a Fourier-spectral method, Phys. D.  179 (3-4) (2003)  211-228.

 \bibitem{C. Liu(2017)}
C.Liu, and H. Wu,  An energetic variational approach for the Cahn-Hilliard equation with dynamic boundary conditions, arXiv preprint arXiv:1710.08318, 2017.

\bibitem{M. Mimura(1984)}
M. Mimura, T. Nakaki,  K. Tomoeda, A numerical approach to interface curves for some nonlinear diffusion equations, Japan J. Appl. Math. 1 (1984) 93-139.

\bibitem{T. Nakaki(2003)}
T. Nakaki, K. Tomoeda,  Numerical approach to the waiting time for the one-dimensional porous medium equation, Quart. Appl. Math.  61 (4) (2002) 601-612.


 \bibitem{Y. Nesterov(1994)}
Y. Nesterov, A. Nemirovskii,  Interior-point polynomial algorithms in convex programming, SIAM, 1994. 

\bibitem{C. Ngo(2017)}
C. Ngo, W.Z.  Huang,  A study on moving mesh finite element solution of the porous medium equation,  J. Compu. Phys. 331 (2017) 357-380.


 \bibitem{O.A. Oleinik(1958)}
O.A. Ole\v{\i}nik, A.S. Kala\v{s}inkov,  Y. \v{C}\v{z}ou, The Cauchy problem and boundary problems for equations of the type of non-stationary filtration, Izv. Akad. Nauk SSSR, Ser. Mat. 22 (1958) 667-704.

\bibitem{L. Onsager(1931)}
L. Onsager, Reciprocal relations in irreversible processes, Phys. Rev., II. Ser. 38 (1931) 2265-2279.

\bibitem{L. Onsager1(1931)}
L. Onsager, Reciprocal relations in irreversible processes, Phys. Rev., I. 37 (4) (1931) 405.

\bibitem{R.E. Pattle (1959)}
R.E. Pattle,  Diffusion from an instantaneous point source with concentration dependent coefficient, Quart. J. Mech. Appl. Math. 12 (1959)  407-409.


%


\bibitem{S. Shmarev(2003)}
S.I. Shmarev, Interfaces in multidimensional diffusion equations with absorption terms, Nonlinear Anal. 53 (2003) 791-828.


\bibitem{S. Shmarev(2005)}
S. Shmarev, Interfaces in solutions of diffusion-absorption equations in arbitrary space dimension, in: Trends in Partial Differential Equations of Mathematical Physics, in: Progr. Nonlinear Differential Equations Appl.  Birkh\"auser, Basel, 2005, pp. 257-273.

\bibitem{J. W. Strutt(1873)}
J.W. Strutt, Some general theorems relating to vibrations, P. Lond. Math. Soc. IV (1873) 357-368.

\bibitem{K. Tomoeda(1983)}
K. Tomoeda, M. Mimura, Numerical approximations to interface curves for a porous medium equation, Hiroshima Math. J. 13 (1983) 273-294.

\bibitem{J. L. Vazquez(2007)}
J.L. V\'azquez, The Porous Medium Equation, Oxford University Press, Oxford, 2007.

\bibitem{C. Wang(2000)}
C. Wang, J.-G. Liu, Convergence of gauge method for incompressible flow, Math. Comp.  69 (2000) 1385-1407.


 \bibitem{M. Westdickenberg(2010)}
M. Westdickenberg, J. Wilkening, Variational particle schemes for the porous medium equation and for the system of isentropic Euler equations, ESAIM: M2AN. 44 (1) (2010) 133-166.

\bibitem{Ya.B. Zeldovich (1950)}
Ya.B. Zel¡¯dovich, A.S. Kompaneets, Towards a theory of heat conduction with thermal conductivity depending on the temperature. In Collection of Papers Dedicated to 70th Anniversary of A. F. Ioffe. Izd.
Akad. Nauk SSSR, Moscow  (1950)  61-72.

\bibitem{Q. Zhang(2009)}
Q. Zhang,  Z.L. Wu, Numerical simulation for porous medium equation by local discontinuous Galerkin finite element method, J. Sci. Comput. 38 (2) (2009) 127-148.




\end{thebibliography}

}

\end{document}